# EXTENSIONS OF NUMBER FIELDS WITH WILD RAMIFICATION OF BOUNDED DEPTH

FARSHID HAJIR AND CHRISTIAN MAIRE

ABSTRACT. We consider $p$-extensions of number fields such that the filtration of the Galois group by higher ramification groups is of prescribed finite length. We extend well-known properties of tame extensions to this more general setting; for instance, we show that these towers, when infinite, are "asymptotically good" (an explicit bound for the root discriminant is given). We study the difficult problem of bounding the relation-rank of the Galois groups in question. Results of Gordeev and Wingberg imply that the relation-rank can tend to infinity when the set of ramified primes is fixed but the length of the ramification filtration becomes large. We show that all $p$-adic representations of these Galois groups are potentially semistable; thus, a conjecture of Fontaine and Mazur on the structure of tamely ramified Galois $p$-extensions extends to our case. Further evidence in support of this conjecture is presented.

## Introduction

Fix a prime number $p$, a number field $K$, and a finite set $S$ of primes of $K$. Let $S_p$ be the set of all primes of $K$ of residue characteristic $p$. Inside a fixed algebraic closure $\overline{K}$ of $K$, let $K_S$ be the maximal $p$-extension (Galois extension with pro-$p$ Galois group) of $K$ unramified outside $S$, and put $G_S = \text{Gal}(K_S/K)$. The study of these "fundamental groups" is governed by a dichotomy between the tame ($S \cap S_p = \emptyset$) and wild ($S \cap S_p \neq \emptyset$) cases.

One feature of this dichotomy is the following. In the tame case, every open subgroup of $G_S$ has finite abelianization (following Lubotzky, we say $G_S$ is *FAb*). On the other hand, if $S_p \subseteq S$, then $G_S$ has a surjection onto $\mathbb{Z}_p^{r_2+1}$ (induced by the $\mathbb{Z}_p$-extensions of $K$), where $r_2$ is the number of imaginary places of $K$. (For surjections of $G_S$ to $\mathbb{Z}_p$ when $S \subset S_p$, see [M]). Indeed, the difference between the tame and wild cases is highlighted by a conjecture of Fontaine and Mazur [FM] which predicts that, in the tame case, $G_S$ is "$p$-adically finite," meaning it has no infinite $p$-adic analytic quotients.

A second, and subtly related, feature is the following: for $\mathfrak{p} \in S - S_p$, the filtration $D^{-1}(K_S/K, \mathfrak{p}) \supseteq D^0(K_S/K, \mathfrak{p}) \supseteq \cdots$ of $G_S$ by higher ramification groups at $\mathfrak{p}$ (in the upper numbering) has length at most 2, i.e. $D^1(K_S/K, \mathfrak{p})$ vanishes, whereas in the case of wild ramification in an infinite $p$-extension, it is often the case that the higher ramification groups of all indices are non-trivial; the latter condition is called "deeply ramified," [CG], the archetypal example being a $\mathbb{Z}_p$-extension.

In this paper, we study a generalization of tame extensions, namely towers with wild ramification of bounded "depth." To be precise, let $\nu : S \to [0, \infty]$, sending $\mathfrak{p}$ to $\nu_\mathfrak{p}$, be an arbitrary map (which will serve to limit the depth of ramification). Now define the group $G_{S,\nu}$ as the quotient of $G_S$ by the closed normal subgroup generated by all higher ramification groups $D^{\nu_\mathfrak{p}}(K_S/K, \mathfrak{p})$ as $\mathfrak{p}$ runs over $S$. The fixed field $K_{S,\nu}$ of this subgroup, with Galois group $\operatorname{Gal}(K_{S,\nu}/K) = G_{S,\nu}$, is the compositum of all finite $p$-extensions of $K$ having vanishing $\nu_\mathfrak{p}{}^{th}$ higher ramification group at all $\mathfrak{p} \in S$. If the image of $\nu$ is simply $\nu(S) = \{\infty\}$, then $G_{S,\nu}$ is nothing but $G_S$, and much – most notably an estimate for its relation-rank – is known about this group [Sh], [K2], [NSW], [Gr].

The case of greatest interest for us is when $\nu$ is "finite", i.e. $\nu(S) \subseteq [0, \infty)$. Our philosophy is that "everything" which is true about $G_S$ under the tame condition ($S \cap S_p = \emptyset$) is also true about $G_{S,\nu}$ under the finiteness of depth condition (i.e. $\nu$ is finite). For example, when $\nu$ is finite, $G_{S,\nu}$ is FAb. Moreover, when $\nu$ is finite and $G_{S,\nu}$ is infinite, the number fields in the tower $K_{S,\nu}/K$ form an "asymptotically good" family (in the sense of Tsfasman and Vladut [TV]), i.e. the root discriminant of these number fields remains bounded. The basic reason is that the exponent of a prime $\mathfrak{p}$ in the relative discriminant of $L/K$, where $L$ is contained in $K_S$, is a sum of orders of ramification groups and therefore grows linearly in $[L:K]$. An explicit bound is given in section 4; a similar idea for function fields appears in Perret [P], cf. our Remarks 3.6 and 4.3. Another perspective is that the fields in $K_{S,\nu}/K$ have bounded conductor when $\nu$ is finite; see, for example, Shirai [Shi]. For applications of asymptotically good families in the number field and function field settings, see [AM], [HM1], [HM2], [HM3] and especially [TV].

We investigate the group $G_{S,\nu}$ via generators and relations. While calculating its generator-rank $d(G_{S,\nu}) = \dim_{\mathbb{F}_p} H^1(G_{S,\nu}, \mathbb{F}_p)$ reduces to a standard calculation in class field theory, estimating the relation-rank $r(G_{S,\nu}) = \dim_{\mathbb{F}_p} H^2(G_{S,\nu}, \mathbb{F}_p)$ presents an essential technical difficulty we have not been able to overcome. This is explained in section 5.

We are, however, able to make certain interesting observations regarding the relation-rank of $G_{S,\nu}$. First, there are non-trivial cases (i.e. where $\nu_\mathfrak{p} \geq 2$ for some $\mathfrak{p} \in S \cap S_p$) where $G_{S,\nu}$ is finite, and, of course, $G_{S,\nu}$ can also be infinite since already $G_\emptyset$ can be infinite. However, it is not known whether, for finite $\nu$, $G_{S,\nu}$ is even finitely presentable (has finite relation-rank), though, according to the philosophy explained above, we suspect this is so. Our main result here is the construction (using results of Wingberg and Gordeev) of a family of examples with *infinite* $\nu$ and $p = 2$ achieving $r(G_{S,\nu}) = \infty$. These examples allow us to show that for fixed $S$ and growing, but finite, $\nu$, the relation-rank can become arbitrarily large.

In the final part of the paper, further questions about the structure of $G_{S,\nu}$ are discussed, especially as regards $p$-adic representations of these groups. Let us first examine the tamely ramified case. Consider an irreducible $p$-adic representation $\rho$ of $\operatorname{Gal}(\overline{K}/K)$ which factors through $G_S$, where $S \cap S_p = \emptyset$. A theorem of Grothendieck ensures that $\rho$ is potentially semistable. Fontaine and Mazur [FM] conjecture that potentially semistable Galois representations unramified outside a finite set of places, such as our $\rho$, must arise from the action of the absolute Galois group of $K$ on a subquotient of the étale cohomology of some algebraic variety over $K$; algebro-geometric considerations then imply that $\rho$ has finite image. This is how Fontaine and Mazur arrive at the prediction that the tame groups $G_S$ are $p$-adically finite (since all finitely generated $p$-adic analytic groups are linear over $\mathbb{Z}_p$).

Our philosophy that the structure of Galois groups with wild ramification of bounded depth mirrors that of the tame case led us to the suspicion that $p$-adic representations which vanish at all higher ramification groups of *finite* depth should be potentially semistable. Moreover, following Fontaine and Mazur, we would also expect that, when $\nu$ is finite, $G_{S,\nu}$ is $p$-adically finite. The relevance of a Theorem of Coates and Greenberg (based on an important result of Sen) for these questions was pointed out to us by A. Schmidt [Sch2], whom we thank. Namely, if $(S,\nu)$ is a finitely indexed set for $K$, and $L/K$ is a subextension of $K_{S,\nu}/K$, with $p$-adic analytic Galois group $\operatorname{Gal}(L/K)$, then $L/K$ is potentially tamely ramified, i.e. there exists a finite Galois extension $K'/K$ with $K' \subset L$ such that $L/K'$ is tamely ramified.

Three immediate corollaries are: 1. A generalization of Grothendieck's theorem, namely: Every $p$-adic representation of $\operatorname{Gal}(\overline{K}/K)$ factoring through $G_{S,\nu}$ (with $\nu$ finite) is potentially semistable; 2. The Fontaine-Mazur conjecture implies that for finite $\nu$, the groups $G_{S,\nu}$ are $p$-adically finite; and 3. The Fontaine-Mazur conjecture also implies that if $L/K$ is a Galois extension with infinite $p$-adic analytic Galois group $\operatorname{Gal}(L/K)$, then $L/K$ is either ramified at infinitely many primes, or it is deeply ramified at some prime $\mathfrak{p}$ of residue characteristic $p$, meaning the ramification groups of all indices are non-trivial. Examples of extensions with $p$-adic analytic Galois group ramified at infinitely many primes were recently constructed by Ramakrishna [R].

It is worth noting that in case $L/K$ is ramified at infinitely many primes, a simple calculation yields that the root discriminants of the fields in this extension tend to infinity. It is natural to ask, then, if the same is true when $L/K$ is infinitely ramified in the other sense, i.e. is deeply ramified. We answer this question in the affirmative in Section 4. Thanks to the results described in the preceding paragraph and the fact that tame extensions are asymptotically good, one can then reformulate the Fontaine-Mazur conjecture as follows: If $K$ is a number field and $L/K$ is a Galois extension such that $\operatorname{Gal}(L/K)$ is an infinite $p$-adic analytic group, then $L/K$ is asymptotically bad. In this way, one can interpret the conjecture as a statement about the growth rate of the index of the $n$th higher ramification groups inside deeply ramified $p$-adic analytic $p$-extensions. A more far-reaching related question suggested by the extension of Grothendieck's semistability theorem to the finite-depth case is to investigate the possibility of characterizing semistable $p$-adic representations in terms of the growth rate of the index of the corresponding higher ramification groups.

In the final section, we give some examples where we are able to check that $G_{S,\nu}$ is $p$-adically finite without having to assume the tamely ramified Fontaine-Mazur conjecture. (But we are not able to check whether these $G_{S,\nu}$ are actually infinite!)

Since much of our discussion holds in the context of all algebraic extensions, not just $p$-extensions, we should explain that we have restricted ourselves to $p$-extensions here partly to fix ideas, but also partly because, for number fields, even the tame situation is not well understood in the more general setting. For example, it is not known whether the Galois group over $K$ of the maximal algebraic extension of $K$ unramified outside a finite set $S$ is finitely generated or not.

The organization of the paper is as follows. The first two sections comprise a preliminary chapter of definitions and properties of ramification groups. We define the extensions of interest to us in section 3, and calculate the behavior of root discriminants in them in section 4. In section 5, we use theorems of Gordeev and Wingberg to explore the structure of the groups $G_{S,\nu}$. In particular, we pinpoint a difficulty in directly extending the method for calculating the relation rank of $G_S$ to that of $G_{S,\nu}$, and we construct a family of $G_{S,\nu}$ with infinite $\nu$ which are not finitely presentable; this furnishes examples of growing but finite $\nu$, for a fixed $S$, for which the relation rank (and partial Euler characteristic $r - d$) of $G_{S,\nu}$ tend to infinity. Sections 6 and 7 are concerned with $p$-adic representations of $G_{S,\nu}$; in particular, we state an extension of Conjecture 5a of [FM], for which we give some supporting evidence (Corollaries 7.9 and 7.10).

We would like to thank Wayne Aitken, Brian Conrad, Ivan Fesenko, Nick Katz, Franz Lemmermeyer, Barry Mazur, Siman Wong, and especially Alexander Schmidt for helpful conversations and correspondence. The first author was partially supported by a CSUSM Faculty Development Grant. The second author was supported by a Postdoctoral fellowship at MSRI [1] in Berkeley, by CNRS, and by Laboratoire A2X (Bordeaux I).

**Part** 1. **Preliminaries**

## 1. THE SETUP

Let $p$ be a prime number. Consider a number field $K$, a finite set $S$ of primes of $K$ equipped with an indexing function $\nu : S \to [0, \infty]$ sending $\mathfrak{p}$ to $\nu_{\mathfrak{p}}$; we will call the pair $(S, \nu)$ an *indexed set* for $K$. If $\nu(S) \subseteq [0, \infty)$, we will call $(S, \nu)$ a *finitely indexed set*, or simply say that $\nu$ is finite. Sometimes it will be convenient to extend $\nu$ to all places of $K$ by setting $\nu_{\mathfrak{p}} = 0$ for $\mathfrak{p} \notin S$. If $\nu_{\mathfrak{p}} = \infty$ for all $\mathfrak{p} \in S \cap S_p$, then the pair $(S, \nu)$ can be abbreviated by $S$.

We fix once and for all an algebraic closure $\overline{K}$; an "extension of $K$" will mean a subfield of this fixed algebraic closure. Let $\mathfrak{p}$ be a prime ideal of $K$. We will use the following notation.

- $S_p$ is the set of primes of $K$ of residue characteristic $p$; $S - S_p$ means $S - (S \cap S_p)$.
- $K_{\mathfrak{p}}$ is the completion of $K$ at $\mathfrak{p}$;
- $\mathcal{O}_{\mathfrak{p}}$ the ring of integers of $K_{\mathfrak{p}}$;
- $U_{\mathfrak{p}}$ the group of units of $\mathcal{O}_{\mathfrak{p}}$;
- $U_{\mathfrak{p}}^{(i)} = \{u \in U_{\mathfrak{p}} | u - 1 \in \mathfrak{p}^{\lceil i \rceil}\}$ for $i \geq 0$; $U_{\mathfrak{p}}^{(0)} = U_{\mathfrak{p}}$; $U_{\mathfrak{p}}^{(\infty)} = \{1\}$;
- $J_K$ is the group of idèles of $K$;
- $\mathcal{U}_S = \prod_{\mathfrak{p} \notin S} U_{\mathfrak{p}}$; $\mathcal{U}_{S,\nu} = \prod_{\mathfrak{p}} U_{\mathfrak{p}}^{(\nu_{\mathfrak{p}})}$; these products are extended over the infinite places $\mathfrak{p}$ as well, where, for such a place, $U_{\mathfrak{p}} = U_{\mathfrak{p}}^{(0)} = K_{\mathfrak{p}}^{\times}$.
- $\hat{K}$ is the maximal $p$-extension of $K$, with Galois group $\mathcal{G} = \mathrm{Gal}(\hat{K}/K)$;
- $\hat{K}_{\mathfrak{p}}$ is the maximal $p$-extension of $K_{\mathfrak{p}}$ in a fixed Galois closure $\overline{K_{\mathfrak{p}}}$ of $K_{\mathfrak{p}}$, with Galois group $D_{\mathfrak{p}} = \mathrm{Gal}(\hat{K}_{\mathfrak{p}}/K_{\mathfrak{p}})$;
- $N_{S,\nu}$ is the closed normal subgroup of $\mathcal{G}$ generated by all higher ramification groups $D^{\nu_{\mathfrak{p}}}(\hat{K}/K, \mathfrak{p})$ for all $K$-primes $\mathfrak{p}$;
- $G_{S,\nu} = \mathcal{G}/N_{S,\nu} = \mathrm{Gal}(K_{S,\nu}/K)$ is the Galois group of the maximal $p$-extension of $K$ unramified outside $S$ and with ramification of depth at most $\nu_{\mathfrak{p}}$ for $\mathfrak{p} \in S$;
- $d(G) = p\text{-rk}G = \dim_{\mathbb{F}_p} G/[G,G]G^p$ is the $p$-rank of a group $G$;
- for a pro-$p$ group, $r(G) = \dim_{\mathbb{F}_p} H^2(G, \mathbb{F}_p)$ is the minimal number of relations for defining $G$ as a pro-$p$ group;
- $d_{S,\nu} = d(G_{S,\nu}), r_{S,\nu} = r(G_{S,\nu})$;

---

[1]Research at MSRI is supported in part by NSF grant DMS-9701755

- $\zeta_p$ is a primitive $p$th root of unity;
- $\delta_{\mathfrak{p}} = 1$ if $K_{\mathfrak{p}}$ contains $\zeta_p$, 0 otherwise;
- $\delta(\mathfrak{p}, \nu_{\mathfrak{p}})$ is 1 if $\zeta_{\mathfrak{p}} \in U_{\mathfrak{p}}^{(\nu_{\mathfrak{p}})}$, and is 0 otherwise;
- $\delta_K = 1$ if $K$ contains $\zeta_p$, 0 otherwise;
- $\theta_S = 0$ unless $\delta_K = 1$ and $S$ is empty, in which case, $\theta_S = 1$.

Let $S'$ be another finite set of prime ideals of $K$. For later use, we introduce the groups

- $\Delta = \{x \in K^{\times} |\ (x)$ is a $p$th power in the group of fractional ideals of $K\}$
- $\Delta_S = \{x \in \Delta |\ x \in K_{\mathfrak{p}}^{\times p}\ \forall \mathfrak{p} \in S\}/K^{\times p}$;
- $\Delta_{S,\nu} = \{x \in \Delta |\ x \in K_{\mathfrak{p}}^{\times p} U_{\mathfrak{p}}^{(\nu_{\mathfrak{p}})}\ \forall \mathfrak{p} \in S\}/K^{\times p}$;
- $\Delta_S^{S'} = \{x \in K^{\times} |\ x \in U_{\mathfrak{p}} K_{\mathfrak{p}}^{\times p}\ \forall$ finite primes $\mathfrak{p} \notin S'$ and $x \in K_{\mathfrak{p}}^{\times p}\ \forall \mathfrak{p} \in S\}/K^{\times p}$.

Using Kummer theory, one has

**Proposition 1.1.** *Suppose $\zeta_p \in K$. Then,*

(1) *If $S_p \subseteq S$, then $p$-$\mathrm{rk}(\Delta_S) = p$-$\mathrm{rk}(\mathrm{Cl}_K^{S,+})$, where $\mathrm{Cl}_K^{S,+}$ is the $S$-class group of $K$ in the narrow sense.*
(2) *If $S' \cup S = S_p$, then $p$-$\mathrm{rk}(\Delta_S^{S'}) = p$-$\mathrm{rk}(G_{S'}^S)$, where $G_{S'}^S$ is the Galois group over $K$ of the maximal $p$-extension of $K$, unramified outside $S'$ in the narrow sense, in which $S$ splits completely.*

## 2. Higher Ramification groups

For the convenience of the reader, we collect here some definitions and standard properties of higher ramification groups. The reader is invited to consult this section only as needed, and is referred to Serre [Se1] for more details.

### 2.1. Ramification groups with lower numbering.
Consider a finite Galois extension $L/K$ of number fields and a prime $\mathfrak{P}$ of $L$ lying over $\mathfrak{p}$ of $K$.

**Definition 2.1.** The decomposition group $D(L/K, \mathfrak{P}) \subseteq \mathrm{Gal}(L/K)$ of $L/K$ at $\mathfrak{P}$ is the stabilizer of $\mathfrak{P}$ and can be identified with the local Galois group $D(\mathfrak{P}/\mathfrak{p}) = \mathrm{Gal}(L_{\mathfrak{P}}/K_{\mathfrak{p}})$. For $x \geq -1$, we define the higher ramification group of index $x$ in the lower numbering $D_x(L/K, \mathfrak{P}) \subseteq D(\mathfrak{P}/\mathfrak{p})$ by

$$D_x(L/K, \mathfrak{P}) = \{\sigma \in D(\mathfrak{P}/\mathfrak{p}) |\ v_{\mathfrak{P}}(\sigma(\alpha) - \alpha) \geq x + 1, \forall \alpha \in \mathcal{O}_L\}.$$

In the sequel, we will sometimes abbreviate the data $(L/K, \mathfrak{P})$ as $(\mathfrak{P}/\mathfrak{p})$, for instance $D_x(\mathfrak{P}/\mathfrak{p})$ in place of $D_x(L/K, \mathfrak{P})$, etc. Note that the ramification groups of index $-1, 0$ are the decomposition and inertia groups of $L/K$ at $\mathfrak{P}$, respectively. The higher ramification groups give a finite decreasing filtration of the decomposition group: $D_x(L/K, \mathfrak{P}) \subseteq D_y(L/K, \mathfrak{P})$, for $x \geq y$ and $D_x(L/K, \mathfrak{P}) = \{1\}$ for all large enough $x$.

### 2.2. Ramification groups with upper numbering.
Consider the map $\varphi = \varphi_{L/K, \mathfrak{P}} = \varphi_{\mathfrak{P}/\mathfrak{p}}$ from $[-1, \infty)$ to $[-1, \infty)$ defined as follows: For $-1 \leq u \leq 0$, $\varphi(u) = u$; for $u \geq 0$, let $m = \lfloor u \rfloor$ and put

$$\varphi(u) = \frac{1}{g_0} \left(g_1 + \cdots + g_m + (u - m) g_{m+1}\right),$$

where $g_i = |D_i(L/K, \mathfrak{P})|$; in particular, for integral $u$, we have

$$\varphi(u) = -1 + \frac{1}{g_0} \sum_{i=0}^{u} g_i.$$

**Proposition 2.2.** *The function $\varphi$ is continuous, piecewise linear, strictly increasing, concave, and satisfies $\varphi(0) = 0, \varphi(x) \leq x$.*

Let $\psi = \psi_{L/K,\mathfrak{P}} = \psi_{\mathfrak{P}/\mathfrak{p}} : [-1, \infty) \to [-1, \infty)$ be the inverse of $\varphi$.

**Proposition 2.3.** *The map $\psi$ is continuous, piecewise linear, strictly increasing, convex, and satisfies $\psi(0) = 0, \psi(y) \geq y$. If $n$ is an integer, then $\psi(n)$ is an integer.*

**Definition 2.4.** We define the ramification group of index $y \geq -1$ in the upper numbering by $D^y(L/K, \mathfrak{P}) = D_{\psi(y)}(L/K, \mathfrak{P})$.

*Remark* 2.5. We note that $D_x(L/K, \mathfrak{P}) = D^{\varphi(x)}(L/K, \mathfrak{P})$, $D^y(L/K, \mathfrak{P}) \subseteq D_y(L/K, \mathfrak{P})$, and for $y \geq z$, $D^y(L/K, \mathfrak{P}) \subseteq D^z(L/K, \mathfrak{P})$. The ramification groups of $\mathfrak{P}^\sigma$, $\sigma \in \text{Gal}(L/K)$, are the conjugates by $\sigma$ of the ramification groups for $\mathfrak{P}$. We let $D_x(L/K, \mathfrak{p})$ ($D^y(L/K, \mathfrak{p})$) denote the conjugacy class in $\text{Gal}(L/K)$ of $D_x(L/K, \mathfrak{P})$ ($D^y(L/K, \mathfrak{P})$).

Now let us consider the behavior of these groups in a tower $L'/L/K$ where $L'/K$ and $L/K$ are finite Galois extensions, $\mathfrak{P}'$ is a prime of $L'$ lying over $\mathfrak{P}$ of $L$ and $\mathfrak{p}$ of $K$.

**Proposition 2.6.** *We have*

(1) $\varphi_{\mathfrak{P}'/\mathfrak{p}} = \varphi_{\mathfrak{P}/\mathfrak{p}} \circ \varphi_{\mathfrak{P}'/\mathfrak{P}}$, *and* $\psi_{\mathfrak{P}'/\mathfrak{p}} = \psi_{\mathfrak{P}'/\mathfrak{P}} \circ \psi_{\mathfrak{P}/\mathfrak{p}}$;
(2) $D_x(\mathfrak{P}'/\mathfrak{P}) = D_x(\mathfrak{P}'/\mathfrak{p}) \cap \text{Gal}(L'/L)$;
(3) $D^y(\mathfrak{P}/\mathfrak{p}) = D^y(\mathfrak{P}'/\mathfrak{p})D(\mathfrak{P}'/\mathfrak{P})/D(\mathfrak{P}'/\mathfrak{P})$; *in other words, $D^y(\mathfrak{P}/\mathfrak{p})$ is the image of $D^y(\mathfrak{P}'/\mathfrak{p})$ in $\text{Gal}(L/K)$ under the restriction map.*

The groups in the lower numbering behave well under taking subgroups, whereas the groups in the upper numbering behave well under taking quotients. As noted in [Se1], the restriction property allows one to define the ramification groups $D^y(L/K, \mathfrak{P})$ with upper numbering for all profinite extensions $L/K$. We now list some simple consequences of the above proposition, which we will use later.

**Proposition 2.7.**  (1) *If $\mathfrak{P}/\mathfrak{p}$ is unramified, then $D^y(\mathfrak{P}/\mathfrak{p}) = D^y(\mathfrak{P}'/\mathfrak{P})$.*
(2) *If $D^y(\mathfrak{P}'/\mathfrak{P})$ and $D^y(\mathfrak{P}/\mathfrak{p})$ are trivial, then $D^y(\mathfrak{P}'/\mathfrak{p})$ is also trivial.*
(3) *If $L''$ is a Galois extension of $K$ such that $LL'' = L'$, and $D^y(\mathfrak{P}/\mathfrak{p}) = D^y(\mathfrak{P}''/\mathfrak{p}) = \{1\}$, where $\mathfrak{P}'' = \mathfrak{P}' \cap L''$, then $D^y(\mathfrak{P}'/\mathfrak{p}) = \{1\}$.*

*Proof.* 1. As $\mathfrak{P}/\mathfrak{p}$ is unramified, the lower numbering ramification groups of $\mathfrak{P}'/\mathfrak{P}$ are the same as those of $\mathfrak{P}'/\mathfrak{p}$. The map $\psi_{\mathfrak{P}'/\mathfrak{p}}$ is then the same as $\psi_{\mathfrak{P}'/\mathfrak{P}}$. Then one has:

$$\begin{aligned} D^y(\mathfrak{P}'/\mathfrak{p}) &= D_{\psi_{\mathfrak{P}'/\mathfrak{p}}(y)}(\mathfrak{P}'/\mathfrak{p}) = D_{\psi_{\mathfrak{P}'/\mathfrak{P}}(y)}(\mathfrak{P}'/\mathfrak{p}) \\ &= D_{\psi_{\mathfrak{P}'/\mathfrak{P}}(y)}(\mathfrak{P}'/\mathfrak{P}) = D^y(\mathfrak{P}'/\mathfrak{P}). \end{aligned}$$

2. The restriction property shows that $D^y(\mathfrak{P}'/\mathfrak{p}) \subseteq \text{Gal}(L'/L)$. Moreover $D^y(\mathfrak{P}'/\mathfrak{p}) = D_{\psi_{\mathfrak{P}'/\mathfrak{p}}(y)}(\mathfrak{P}'/\mathfrak{p})$. Now using the transitivity of $\psi$ and other elementary properties, we have $D_{\psi_{\mathfrak{P}'/\mathfrak{p}}(y)}(\mathfrak{P}'/\mathfrak{p}) \subseteq D_{\psi_{\mathfrak{P}'/\mathfrak{P}}(y)}(\mathfrak{P}'/\mathfrak{p})$. Thus,

$$D_{\psi_{\mathfrak{P}'/\mathfrak{p}}(y)}(\mathfrak{P}'/\mathfrak{p}) \subseteq \text{Gal}(L'/L) \cap D_{\psi_{\mathfrak{P}'/\mathfrak{P}}(y)}(\mathfrak{P}'/\mathfrak{p}) = D_{\psi_{\mathfrak{P}'/\mathfrak{P}}(y)}(\mathfrak{P}'/\mathfrak{P}) = D^y(\mathfrak{P}'/\mathfrak{P}).$$

The triviality of $D^y(\mathfrak{P}'/\mathfrak{P})$ allows us to conclude.

3. Follows easily from the restriction property.  □

**Proposition 2.8.** *Suppose $L/K$ is a finite $p$-extension.*

*1) If $\mathfrak{p}$ has residue characteristic $\ell \neq p$, $D^y(\mathfrak{P}/\mathfrak{p})$ vanishes for all $y > 0$.*

*2) If $\mathfrak{p}$ has residue characteristic $p$, and $D^y(\mathfrak{P}/\mathfrak{p})$ vanishes for some $y \leq 1$, then $\mathfrak{P}/\mathfrak{p}$ is unramified.*

*Proof.* 1) Since $x = \psi_{\mathfrak{P}/\mathfrak{p}}(y) > 0$, we have $D^y(\mathfrak{P}/\mathfrak{p}) = D_x(\mathfrak{P}/\mathfrak{p}) = D_1(\mathfrak{P}/\mathfrak{p})$. Since $L/K$ is a $p$-extension, $\mathfrak{P}/\mathfrak{p}$ is at most tamely ramified and the wild inertia group $D_1(\mathfrak{P}/\mathfrak{p})$ is trivial.

2) We claim that $x := \psi_{\mathfrak{P}/\mathfrak{p}}(1) = 1$. We already know that $x$ is a positive integer. Let us show that $x \leq 1$. We have
$$1 = \varphi_{\mathfrak{P}/\mathfrak{p}}(x) = \frac{g_1 + \cdots + g_x}{g_0},$$
where $g_i = |D_i(\mathfrak{P}/\mathfrak{p})|$. As $D_0(\mathfrak{P}/\mathfrak{p})/D_1(\mathfrak{P}/\mathfrak{p})$ is trivial, $g_0 = g_1$. One obtains:
$$1 + \frac{x-1}{g_0} \leq 1,$$
and so $x = 1$. Now we see that if $y \leq 1$, then
$$D_1(\mathfrak{P}/\mathfrak{p}) = D_{\psi_{\mathfrak{P}/\mathfrak{p}}(1)}(\mathfrak{P}/\mathfrak{p}) \subseteq D_{\psi_{\mathfrak{P}/\mathfrak{p}}(y)}(\mathfrak{P}/\mathfrak{p}) = D^y(\mathfrak{P}/\mathfrak{p}) = \{1\},$$
and this implies the triviality of $D_0(\mathfrak{P}/\mathfrak{p})$ because this is a $p$-group and $D_0(\mathfrak{P}/\mathfrak{p})/D_1(\mathfrak{P}/\mathfrak{p})$ has order prime to $p$.   $\square$

*Remark* 2.9. If $\mathfrak{P}$ has residue characteristic $p$, $\psi_{\mathfrak{P}/\mathfrak{p}}(y) = y$ for $y \in [0, 1]$.

**Proposition 2.10.** *Suppose $L/K$ is an abelian extension. Then $D_1(\mathfrak{P}/\mathfrak{p})$ vanishes if and only if $D^1(\mathfrak{P}/\mathfrak{p})$ does.*

*Proof.* One direction is trivial because $D^1(\mathfrak{P}/\mathfrak{p}) \subseteq D_1(\mathfrak{P}/\mathfrak{p})$. We know that the image of principal units under the local reciprocity map generates $D^1(\mathfrak{P}/\mathfrak{p})$. Thus, if this group is trivial, by using the structure of the group of local units we deduce the fact that the inertia group has order prime to the residue characteristic of $\mathfrak{p}$, and then $D_1(\mathfrak{P}/\mathfrak{p})$ is trivial.   $\square$

**Proposition 2.11.** *Suppose $D^y(\mathfrak{P}'/\mathfrak{p}) = \{1\}$ for some $y$. Then $D^z(\mathfrak{P}'/\mathfrak{P}) = \{1\}$ for $z = \psi_{\mathfrak{P}/\mathfrak{p}}(y)$.*

*Proof.* One has:
$$\begin{aligned} D^z(\mathfrak{P}'/\mathfrak{P}) &= D_{\psi_{\mathfrak{P}'/\mathfrak{P}}(z)}(\mathfrak{P}'/\mathfrak{P}) = D_{\psi_{\mathfrak{P}'/\mathfrak{p}}(y)}(\mathfrak{P}'/\mathfrak{P}) \\ &= D_{\psi_{\mathfrak{P}'/\mathfrak{p}}(y)}(\mathfrak{P}'/\mathfrak{p}) \cap \mathrm{Gal}(L'/L) \\ &= D^y(\mathfrak{P}'/\mathfrak{p}) \cap \mathrm{Gal}(L'/L) \\ &= \{1\}. \end{aligned}$$
$\square$

## Part 2. **Wild ramification of bounded depth**

### 3. Definition and simple properties of towers with bounded ramification

Let $S$ be a finite set of primes of $K$, equipped with a map $\nu : S \to [0, \infty]$ sending $\mathfrak{p}$ to $\nu_{\mathfrak{p}}$.

**Definition 3.1.** Suppose $L/K$ is a finite Galois extension and $\mathfrak{p}$ is a prime of $K$. We say that the ramification of $L/K$ at $\mathfrak{p}$ is of depth at most $y$ if $D^y(\mathfrak{P}/\mathfrak{p})$ vanishes for all primes $\mathfrak{P}$ of $L$ lying over $\mathfrak{p}$.

The ramification of $L/K$ at $\mathfrak{p}$ is of depth at most 0 (at most 1) means that $\mathfrak{p}$ is unramified (at most tamely ramified) in $L/K$.

**Definition 3.2.** Let $K_{S,\nu}$ be the compositum, inside our fixed algebraic closure of $K$, of all finite $p$-extensions $L$ of $K$ which are unramified outside $S$ and have the further property that, for every $\mathfrak{p} \in S$, the ramification of $L/K$ at $\mathfrak{p}$ is of depth at most $\nu_\mathfrak{p}$. Note that the infinite places of $K$ are required to split completely in $K_{S,\nu}$. Let $G_{S,\nu} = \mathrm{Gal}(K_{S,\nu}/K)$.

By Proposition 2.7, $K_{S,\nu}$ is the maximal $p$-extension of $K$ unramified outside $S$ and ramified to depth at most $\nu_\mathfrak{p}$ for every $\mathfrak{p} \in S$. If $K_S$ is the maximal $p$-extension of $K$ unramified outside $S$, then $K_{S,\nu}$ is the fixed field corresponding to the closed normal subgroup of $G_S = \mathrm{Gal}(K_S/K)$ generated by all $D^{\nu_\mathfrak{p}}(K_S/K,\mathfrak{p})$. If $\nu$ vanishes identically, then $K_{S,\nu}$ is simply the Hilbert $p$-class field tower of $K$. By Proposition 2.8, we may assume without loss of generality that

$$\mathfrak{p} \in S - S_p \Rightarrow \nu_\mathfrak{p} = 1, \qquad \mathfrak{p} \in S \cap S_p \Rightarrow \nu_\mathfrak{p} > 1.$$

*Example* 3.3. Take $K = \mathbb{Q}, S = \{p\}$. Then $G_S \simeq \mathbb{Z}_p$. For any finite indexing $\nu$ of $S$, we have $G_{S,\nu}$ is finite with vanishing Euler characteristic.

**Definition 3.4.** Suppose $L/K$ is a finite extension contained in $K_{S,\nu}$. We lift the indexed set $(S,\nu)$ of $K$ to an indexed set $(S,\nu)_L = (S(L),\nu(L))$ of $L$ as follows:

$$S(L) = \{\mathfrak{P} \subseteq \mathcal{O}_L | \mathfrak{P} \text{ divides } \mathfrak{p}\mathcal{O}_L \text{ for some } \mathfrak{p} \in S\}, \qquad \nu(L)_\mathfrak{P} = \psi_{\mathfrak{P}/\mathfrak{p}}(\nu_\mathfrak{p}).$$

Let us note that for a tower $F/L/K$, the indexed sets $(S,\nu)_F$ and $(S(L),\nu(L))_F$ coincide, thanks to the transitivity of $\psi$.

Suppose $\nu$ is *finite*, and $K_1$ is the maximal abelian extension of $K$ contained in $K_{S,\nu}$; this is the field associated by class field theory to the idèle subgroup:

$$\mathcal{U}_{S,\nu} = \prod_{\mathfrak{p} \notin S} U_\mathfrak{p} \prod_{\mathfrak{p} \in S} U_\mathfrak{p}^{(\nu_\mathfrak{p})}.$$

Thus, $K_1$ is the ray class field of $K$ modulo $\mathfrak{m}_{S,\nu} = \prod_{\mathfrak{p} \in S} \mathfrak{p}^{\lceil \nu_\mathfrak{p} \rceil}$; in particular, $K_1/K$ is finite.

To simplify notation, let us write $(S(1),\nu(1))$ for $(S(K_1),\nu(K_1))$. Define inductively a tower of finite abelian extensions as follows: $K_{n+1}$ is the maximal abelian extension of $K_n$ contained in $K_{S(n),\nu(n)}$ and $(S(n+1),\nu(n+1))$ is the lift of $(S(n),\nu(n))$ to $K_{n+1}$. Let $K_\infty = \cup_n K_n$. By maximality for each step, $K_\infty/K$ is a Galois pro-$p$ extension.

**Theorem 3.5.** *Suppose $(S,\nu)$ is a finitely indexed set for $K$. With notations as above, we have*

(1) $K_\infty = K_{S,\nu}$;
(2) *If $K/K_0$ is a Galois extension such that $(S,\nu)$ is $\mathrm{Gal}(K/K_0)$-stable, then all $K_n/K_0$, and $K_{S,\nu}/K_0$, are Galois extensions;*
(3) *If $L/K$ is a finite Galois extension contained in $K_{S,\nu}$, then $K_{S,\nu} = L_{S(L),\nu(L)}$;*
(4) *The open subgroups of $G_{S,\nu}$ have finite abelianization, i.e. $G_{S,\nu}$ is FAb.*

*Proof.* 1) For $\mathfrak{p} \in S$, choose a compatible system $\mathfrak{P}_n \in S(n)$, lying over $\mathfrak{p}$. Recall that $\nu_n(L)_{\mathfrak{P}_n} = \psi_{\mathfrak{P}_n/\mathfrak{p}}(\nu_\mathfrak{p})$. To increase the readability of what follows, let us use the temporary notation $D[j](\mathfrak{P}/\mathfrak{p}) = D^j(\mathfrak{P}/\mathfrak{p})$ and $D\{j\}(\mathfrak{P}/\mathfrak{p}) = D_j(\mathfrak{P}/\mathfrak{p})$ for the upper and lower ramification groups respectively. Consider a positive integer $n$ such that $K_n \subseteq K_{S,\nu}$ (this is the case for $n = 1$). We calculate:

$$\begin{aligned} D[\nu_\mathfrak{p}](\mathfrak{P}_{n+1}/\mathfrak{p}) \cap \mathrm{Gal}(K_{n+1}/K_n) &= D\{\psi_{\mathfrak{P}_{n+1}/\mathfrak{p}}(\nu_\mathfrak{p})\}(\mathfrak{P}_{n+1}/\mathfrak{p}) \cap \mathrm{Gal}(K_{n+1}/K_n) \\ &= D\{\psi_{\mathfrak{P}_{n+1}/\mathfrak{p}}(\nu_\mathfrak{p})\}(\mathfrak{P}_{n+1}/\mathfrak{P}_n) \\ &= D\{\psi_{\mathfrak{P}_{n+1}/\mathfrak{P}_n}(\psi_{\mathfrak{P}_n/\mathfrak{p}}(\nu_\mathfrak{p}))\}(\mathfrak{P}_{n+1}/\mathfrak{P}_n) \\ &= D[\psi_{\mathfrak{P}_n/\mathfrak{p}}(\nu_\mathfrak{p})](\mathfrak{P}_{n+1}/\mathfrak{P}_n) \\ &= \{1\}. \end{aligned}$$

But, by the restriction property of higher ramification groups, and the fact that $D[\nu_\mathfrak{p}](\mathfrak{P}_n/\mathfrak{p})$ vanishes, $D[\nu_\mathfrak{p}](\mathfrak{P}_{n+1}/\mathfrak{p}) \subseteq \mathrm{Gal}(K_{n+1}/K_n)$, hence is trivial. By induction, $K_\infty \subseteq K_{S,\nu}$.

To show the reverse inclusion, it suffices to show that if $L$ is a finite abelian $p$-extension of $K_n$ (for an arbitrary $n$) which is Galois over $K$ and contained in $K_{S,\nu}$, then $L \subseteq K_{n+1}$. If $\mathfrak{P}$ is a prime of $L$ dividing $\mathfrak{P}_n$, we have:

$$\begin{aligned} D[\nu(L)_{\mathfrak{P}_n}](\mathfrak{P}/\mathfrak{P}_n) &= D[\psi_{\mathfrak{P}_n/\mathfrak{p}}(\nu_\mathfrak{p})](\mathfrak{P}/\mathfrak{P}_n) \\ &= D\{\psi_{\mathfrak{P}/\mathfrak{P}_n} \circ \psi_{\mathfrak{P}_n/\mathfrak{p}}(\nu_\mathfrak{p})\}(\mathfrak{P}/\mathfrak{P}_n) \\ &= D\{\psi_{\mathfrak{P}/\mathfrak{p}}(\nu_\mathfrak{p})\}(\mathfrak{P}/\mathfrak{p}) \cap \mathrm{Gal}(L/K_n) \\ &= D[\nu_\mathfrak{p}](\mathfrak{P}/\mathfrak{p}) \\ &= \{1\}. \end{aligned}$$

By the maximality of $K_{n+1}/K_n$, we then have $L \subseteq K_{n+1}$, completing the proof of 1.

2) is clear.

3) Proposition 2.11 gives $K_{S,\nu} \subseteq L_{S(L),\nu(L)}$. We show the reverse inclusion. Let $(L_n)_n$ be the sequence of maximal abelian extensions giving $\cup_n L_n = L_{S(L),\nu(L)}$. Put $L' = L_n$ for an arbitrary fixed $n$. By 2), $L'/K$ is Galois. One gets:

$$\begin{aligned} \{1\} = D[\nu(L)_\mathfrak{P}](\mathfrak{P}'/\mathfrak{P}) &= D\{\psi_{\mathfrak{P}'/\mathfrak{P}}(\nu(L)_\mathfrak{P})\}(\mathfrak{P}'/\mathfrak{P}) \\ &= D\{\psi_{\mathfrak{P}'/\mathfrak{p}}(\nu_\mathfrak{p})\}(\mathfrak{P}'/\mathfrak{P}) \\ &= D\{\psi_{\mathfrak{P}'/\mathfrak{p}}(\nu_\mathfrak{p})\}(\mathfrak{P}'/\mathfrak{p}) \cap \mathrm{Gal}(L'/L) \\ &= D[\nu_\mathfrak{p}](\mathfrak{P}'/\mathfrak{p}) \cap \mathrm{Gal}(L'/L). \end{aligned}$$

Moreover, $D[\nu_\mathfrak{p}](\mathfrak{P}/\mathfrak{p}) = \{1\}$, and so $D[\nu_\mathfrak{p}](\mathfrak{P}'/\mathfrak{p}) \subseteq \mathrm{Gal}(L'/L)$. In conclusion $D[\nu_\mathfrak{p}](\mathfrak{P}'/\mathfrak{p}) = \{1\}$, and $L_n = L' \subseteq K_{S,\nu}$.

4) If $H$ is an open subgroup of $G_{S,\nu}$, its fixed field $L$ is a finite extension of $K$ and the abelianization of $H = G_{S(L),\nu(L)}$ is isomorphic to a ray class group with finite conductor, hence is finite. $\square$

*Remark* 3.6. 1) The naive lift of $\nu$ from $K$ to $L$ namely $\mathfrak{P} \mapsto \nu_\mathfrak{p}$ (for a prime $\mathfrak{P}$ of $L$ lying over $\mathfrak{p}$ of $K$) would lead to examples where $L_{S(L),\nu(L)} \neq K_{S,\nu}$. Here is a simple example for $p = 2$: let $K = \mathbb{Q}$, $L = \mathbb{Q}(\sqrt{3})$ and $F = K(\sqrt{2})$. We let $\mathfrak{P}, \mathfrak{p}$ be the unique primes of $F, L$ (respectively) dividing 2. One checks that $D^3(\mathfrak{P}/2)$ vanishes but that $D^3(\mathfrak{P}/\mathfrak{p})$ does not.

2) We should remark that the towers (of function fields) constructed in Perret [P] are based on this naive lift of $(S, \nu)$. They are contained in $K_{S,\nu}$ by Proposition 2.7, and, thanks to Hasse-Arf, have a slightly lower root discriminant bound than the tower $K_{S,\nu}/K$, cf. Remark 4.3. The Galois group of Perret's tower, however, being less natural, is probably very difficult to study. For example, Neiderreiter and Xing [NX] have shown that the relation-rank estimate conjectured by Perret is not correct, at least for function fields over $\mathbb{F}_2$.

3) By the above theorem, the fields $K_n$ are simply the fields fixed by the "commutator series" of $G_{S,\nu}$. Also, $K_\emptyset$ = the maximal unramified $p$-extension of $K$, is contained in $K_{S,\nu}$, for an arbitrary $(S,\nu)$.

We conclude this section with a calculation of the generator rank of $G_{S,\nu}$. By the Burnside Basis theorem, this reduces to calculating the $p$-rank of its maximal abelian quotient, which can be expressed, by class field theory, as the $p$-rank of an appropriate ray class group. The reduction to abelian extensions also allows us to assume, without loss of generality, that $\nu$ takes integer values (thanks to Hasse-Arf). We also note that if $\min_{\mathfrak{p}\in S}\nu_\mathfrak{p}$ is large enough, then $p$-rk$(G_{S,\nu}) = p$-rk$(G_S)$.

**Theorem 3.7.** *The generator-rank $d_{S,\nu}$ of $G_{S,\nu}$ is the $p$-rank of the ray class group $\mathrm{Cl}_{K,\mathfrak{m}}$ of conductor $\mathfrak{m} = \prod_{\mathfrak{p}\in S}\mathfrak{p}^{\lceil\nu_\mathfrak{p}\rceil}$. It satisfies*

$$\begin{aligned} d_{S,\nu} &= p\text{-rk}\Delta_{S,\nu} - p\text{-rk}E_K + p\text{-rk}\frac{\mathcal{U}_\emptyset}{\mathcal{U}_{S,\nu}\mathcal{U}_\emptyset^p} \\ &= p\text{-rk}\Delta_{S,\nu} - p\text{-rk}E_K + \sum_{\mathfrak{p}\in S-S_p}\delta_\mathfrak{p} + \sum_{\mathfrak{p}\in S\cap S_p}p\text{-rk}\frac{U_\mathfrak{p}^{(1)}}{U_\mathfrak{p}^{(\nu_\mathfrak{p})}}. \end{aligned}$$

*Proof.* The main observation, as mentioned above, is that the fixed field $K_1$ of (the closure of) the commutator subgroup $[G_{S,\nu}, G_{S,\nu}]$ is the field associated by class field theory to the idèle subgroup:

$$\mathcal{U}_{S,\nu} = \prod_{\mathfrak{p}\notin S} U_\mathfrak{p} \prod_{\mathfrak{p}\in S} U_\mathfrak{p}^{(\nu_\mathfrak{p})}.$$

We recall that we defined certain number groups $\Delta_{S,\nu}$ in section 1. We leave to the reader the verification of the exactness of the sequence:

$$1 \longrightarrow \Delta_{S,\nu} \longrightarrow \Delta_\emptyset \longrightarrow \frac{\mathcal{U}_\emptyset}{\mathcal{U}_\emptyset^p \mathcal{U}_{S,\nu}} \longrightarrow \frac{J_K}{\mathcal{U}_{S,\nu}J_K^p} \longrightarrow \frac{J_K}{K^\times \mathcal{U}_\emptyset J_K^p} \longrightarrow 1.$$

Note that $J_K/K^\times \mathcal{U}_\emptyset J_K^p = \mathrm{Cl}_K/\mathrm{Cl}_K^p$. Moreover, it is easy to establish that $p$-rk$\Delta_\emptyset = p$-rk$\mathrm{Cl}_K + p$-rk$E_K$, using the obvious map from $\Delta_\emptyset$ to the ideal classes killed by $p$. Putting all of this together, we obtain the desired formula. □

**Corollary 3.8.** *The generator rank $d_S$ of $G_S$ satisfies:*

$$d_S = p\text{-rk}\Delta_S - p\text{-rk}E_K + \sum_{\mathfrak{p}\in S-S_p}\delta_\mathfrak{p} + \sum_{\mathfrak{p}\in S\cap S_p}p\text{-rk}U_\mathfrak{p}^{(1)}.$$

*Proof.* For $\mathfrak{p}\in S$, take $\nu_\mathfrak{p} = \infty$. This is the well-known formula of Shafarevich [Sh]. □

## 4. Behavior of the root discriminant

**4.1. Extensions of bounded depth.** We now suppose we have a finitely indexed set $(S,\nu)$ of $K$. We will give an upper bound for the root discriminant of the number fields in the tower $K_{S,\nu}/K$.

**Definition 4.1.** Let $K$ be a number field of degree $n$, and discriminant $\operatorname{disc}_K$. The root discriminant $\operatorname{rd}_K$ of $K$ is defined by
$$\operatorname{rd}_K = |\operatorname{disc}_K|^{1/n}.$$

**Theorem 4.2.** *In the p-extension $K_{S,\nu}/K$, the root discriminant is bounded. More precisely, if $L$ is a finite extension of $K$ contained in $K_{S,\nu}$, then*
$$\operatorname{rd}_L \leq \operatorname{rd}_K \prod_{\mathfrak{p} \in S - S_\mathfrak{p}} (N_{K/\mathbb{Q}}\mathfrak{p})^{1/[K:\mathbb{Q}]} \prod_{\mathfrak{p} \in S \cap S_\mathfrak{p}} (N_{K/\mathbb{Q}}\mathfrak{p})^{(\nu_\mathfrak{p}+1)/[K:\mathbb{Q}]}.$$

*Proof.* It suffices to consider $K \subseteq L \subseteq K_{S,\nu}$ such that $L/K$ is Galois. Fix one of the $r_\mathfrak{p}$ primes $\mathfrak{P}$ of $L$ dividing $\mathfrak{p} \in S$; since $L/K$ is Galois, the ramification index $e_{\mathfrak{P}/\mathfrak{p}} = e_\mathfrak{p}$ and residue degree $f_{\mathfrak{P}/\mathfrak{p}} = f_\mathfrak{p}$ depend only on $\mathfrak{p}$. Let $\varphi = \varphi_{\mathfrak{P}/\mathfrak{p}}$, $\psi = \psi_{\mathfrak{P}/\mathfrak{p}}$.

By the restriction property, $D^{\nu_\mathfrak{p}}(\mathfrak{P}/\mathfrak{p})$ is trivial for all $\mathfrak{p}$ in $S$; moreover, $L/K$ is unramified outside $S$. By definition,
$$D^{\nu_\mathfrak{p}}(\mathfrak{P}/\mathfrak{p}) = D_{\psi(\nu_\mathfrak{p})}(\mathfrak{P}/\mathfrak{p}) = \{1\}.$$
Let $n = \psi(\nu_\mathfrak{p})$ and $m = \lfloor n \rfloor$. The valuation $v_\mathfrak{P}(\mathfrak{d}_{L/K})$ at $\mathfrak{P}$ of the different of $L/K$ is given by Hilbert's formula [Se1]:
$$v_\mathfrak{P}(\mathfrak{d}_{L/K}) = g_0 + g_1 + \cdots + g_{m-1} - m,$$
where $g_j = |D_j(\mathfrak{P}/\mathfrak{p})|$; note that $g_j = 1$ for $j \geq m$. In particular, for $\mathfrak{p} \in S - S_p$ the ramification is tame, so we have
$$v_\mathfrak{P}(\mathfrak{d}_{L/K}) = e_\mathfrak{p} - 1.$$
From the definitions of $\psi, \varphi$, one has:
$$\nu_\mathfrak{p} = \varphi(n) = \frac{g_1 + \cdots + g_m + (n-m)}{g_0}.$$
Moreover,
$$\begin{aligned} v_\mathfrak{P}(\mathfrak{d}_{L/K}) &= g_0 + g_1 + \cdots + g_m - (m+1) \\ &= g_0(\nu_\mathfrak{p} + 1) - (m+1) - (n-m) \\ &= g_0(\nu_\mathfrak{p} + 1) - (n+1) \\ &\leq (e_\mathfrak{p} - 1)(\nu_\mathfrak{p} + 1), \end{aligned}$$
since $g_0 = e_\mathfrak{p}$ and $\nu_\mathfrak{p} = \varphi(n) \leq n$.

Putting the local data together and recalling that $e_\mathfrak{p} f_\mathfrak{p} r_\mathfrak{p} = [L:K]$, we have
$$\begin{aligned} N_{L/\mathbb{Q}}\mathfrak{d}_{L/K} &\leq N_{K/\mathbb{Q}} \prod_{\mathfrak{p} \in S - S_p} \mathfrak{p}^{(e_\mathfrak{p}-1)f_\mathfrak{p} r_\mathfrak{p}} \prod_{\mathfrak{p} \in S \cap S_p} \mathfrak{p}^{(\nu_\mathfrak{p}+1)(e_\mathfrak{p}-1)f_\mathfrak{p} r_\mathfrak{p}} \\ &\leq N_{K/\mathbb{Q}} \prod_{\mathfrak{p} \in S - S_p} \mathfrak{p}^{(1-1/e_\mathfrak{p})[L:K]} \prod_{\mathfrak{p} \in S \cap S_p} \mathfrak{p}^{(\nu_\mathfrak{p}+1)(1-1/e_\mathfrak{p})[L:K]}. \end{aligned}$$
But $[L:K] = [L:\mathbb{Q}]/[K:\mathbb{Q}]$ and
$$\operatorname{rd}_L = \operatorname{rd}_K N_{L/\mathbb{Q}}(\mathfrak{d}_{L/K})^{1/[L:\mathbb{Q}]},$$
so we are done. $\square$

*Remark* 4.3. For the analogue of Perret's tower in the number field case, one has the following root discriminant bound:
$$\operatorname{rd}_L \leq \operatorname{rd}_K \prod_{\mathfrak{p} \in S} (N_{K/\mathbb{Q}}\mathfrak{p})^{\lceil \nu_\mathfrak{p} \rceil / [K:\mathbb{Q}]}.$$

The details, which are essentially in [P], are left to the reader.

## 4.2. Deeply ramified extensions.

**Theorem 4.4.** *Suppose $K$ is a number field and $L/K$ is a deeply ramified Galois extension, i.e. for some prime $\mathfrak{p}$ of $K$, the higher ramification groups at $p$ of all indices are non-trivial. Then $L/K$ is asymptotically bad, i.e. there is a sequence of finite subextensions of $L/K$ with root discriminant tending to infinity.*

*Proof.* Let $\mathfrak{P}$ be a prime of $L$ above $\mathfrak{p}$, a deeply ramified prime. We may choose a sequence of finite normal extensions $F_n/K$ inside $L$ such that (i) $\cup F_n = L$ (so $[F_n : K] \to \infty$), and (ii) for all integers $c$, there is an integer $n$ such that $D^c(F_n/K, \mathfrak{P}_n)$, the upper numbering ramification group of index $c$ at $\mathfrak{P}_n = \mathfrak{P} \cap F_n$ for $\mathrm{Gal}(F_n/K)$, is non-trivial. We will show that the root discriminant of $F_n$ goes to infinity with $n$.

As in the proof of Theorem 4.2, we have

$$\mathrm{rd}_{F_n} = \mathrm{rd}_K \left( \prod_{\mathfrak{q}} \mathbb{N}_{K/\mathbb{Q}} \mathfrak{q}^{\tau(F_n/K,\mathfrak{q})} \right)^{1/[K:\mathbb{Q}]}.$$

Here, the product is over primes $\mathfrak{q}$ of $K$, and $\tau$ is defined by $\tau(F_n/K, \mathfrak{q}) = v_{\mathfrak{Q}}(\mathfrak{d}_{F_n/K})/e_{\mathfrak{Q}/\mathfrak{q}}$, where $\mathfrak{Q}$ is an arbitrary prime of $F_n$ over $\mathfrak{q}$, and $e_{\mathfrak{Q}/\mathfrak{q}}$ is the ramification index of $\mathfrak{Q}/\mathfrak{q}$. It suffices to show, therefore, that $\tau(F_n/K, \mathfrak{p})$ goes to infinity with $n$.

Let $c$ be a positive integer, and choose a corresponding $n$ such that $D^c(F_n/K, \mathfrak{P}_n)$ is non-trivial. Let $z_n$ be the largest integer $j$ such that the lower numbering group $D_j(F_n/K, \mathfrak{P}_n)$ is non-trivial. Recall Hilbert's formula, $v_{\mathfrak{P}_n}(\mathfrak{d}_{F_n/K}) = g_0 + \cdots + g_{z_n} - z_n - 1$ (where $g_i = |D_i(F_n/K, \mathfrak{P}_n)|$), and the definition of $\varphi_{\mathfrak{P}_n/\mathfrak{p}}(z_n) = (g_1 + \cdots g_{z_n})/g_0$ (cf. section 2). These give

$$\tau(F_n/K, \mathfrak{p}) = 1 + \varphi_{\mathfrak{P}_n/\mathfrak{p}}(z_n) - (z_n - 1)/g_0.$$

Recall that we have arranged $D^c(F_n/K, \mathfrak{P}_n) = D_{\psi_{\mathfrak{P}_n/\mathfrak{p}}(c)}(F_n/K, \mathfrak{P}_n)$ to be non-trivial. Since $c$ is an integer, $\psi_{\mathfrak{P}_n/\mathfrak{p}}(c)$ is an integer also (Proposition 2.2). By definition of $z_n$, we have $\psi_{\mathfrak{P}_n/\mathfrak{p}}(c) \leq z_n$. Now we apply $\varphi$ to this inequality to obtain $\varphi_{\mathfrak{P}_n/\mathfrak{p}}(z_n) \geq c$. On the other hand, we have $z_n/g_0 = z_n/e_{\mathfrak{P}_n/\mathfrak{p}} \leq e_{\mathfrak{p}/p}/(p-1)$ (Serre [Se1], Chapter IV, section 2, exercise 3).

Summarizing the argument, we have shown that for all integers $c$, there is an integer $n$ such that $\tau(F_n/K, \mathfrak{p}) \geq c - [K : \mathbb{Q}] - 1$, so $\tau(F_n/K, \mathfrak{p})$ tends to infinity. □

## 5. The relation-rank of $G_S$ and $G_{S,\nu}$

We put

$$r_S = \dim_{\mathbb{F}_p} H^2(G_S, \mathbb{F}_p), \qquad r_{S,\nu} = \dim_{\mathbb{F}_p} H^2(G_{S,\nu}, \mathbb{F}_p),$$

for the relation-rank of $G_S$ and $G_{S,\nu}$, respectively. An estimate for $r_S$ was given by Shafarevich [Sh] in 1963. Estimating $r_{S,\nu}$ appears to be a more delicate problem. In fact, for $p = 2$, using results of Gordeev and Wingberg, we can find examples, with $\nu_{\mathfrak{p}}$ finite for some $\mathfrak{p}$ and infinite for others, where $G_{S,\nu}$ is not finitely presentable! The main problem, then, is the following.

**Question 5.1.** a) Do we have $r_{S,\nu} < \infty$ for finitely indexed sets $(S, \nu)$?

b) If so, give an explicit upper bound $r_{S,\nu} \leq f(S, \nu)$.

Below, we will first indicate briefly why the proof of Shafarevich does not easily generalize to our case. The main results of this section are the construction (for $p = 2$) of examples with $r_{S,\nu} = \infty$ (where $(S,\nu)$ is not finitely indexed), which then leads to a family of examples showing that $r_{S,\nu}$ can tend to infinity for finitely indexed sets $(S,\nu)$, where $S$ is fixed, but $\max_{\mathfrak{p} \in S} \nu_\mathfrak{p}$ tends to infinity.

5.1. **Shafarevich's relation-rank bound.** We saw in the previous section that the generator-rank of $G_{S,\nu}$ is easily recognizable as the $p$-rank of a ray class group. The fundamental arithmetic result about $G_S$ is the following estimate for its relation-rank (rather its partial Euler characteristic $r_S - d_S$) due to Shafarevich [Sh].

**Theorem 5.2.** *We have*
$$r_S - d_S \leq p\text{-rk} E_K - \delta_K + \theta_S - \sum_{\mathfrak{p} \in S \cap S_p} [K_\mathfrak{p} : \mathbb{Q}_p].$$

We note in passing that when $S_p \subseteq S$, we in fact know the partial Euler characteristic exactly: $r_S - d_S = -(r_2 + 1)$ where $r_2$ is the number of imaginary places of $K$ (cf. Haberland [H] and, for $p = 2$, Schmidt [Sch1]).

We will sketch the proof of Theorem 5.2 given by Koch (cf. [K1], [K2]), indicating a difficulty one encounters upon attempting to generalize it by replacing $S$ with $(S,\nu)$.

Recall the group of $(S,\nu)$-unit idèles is
$$\mathcal{U}_{S,\nu} = \prod_\mathfrak{p} U_\mathfrak{p}^{(\nu_\mathfrak{p})},$$
where we have extended $\nu$ to all places by setting $\nu_\mathfrak{p} = 0$ for $\mathfrak{p} \notin S$.

The local reciprocity map $\rho_\mathfrak{p}$ for a local Galois extension $L_\mathfrak{P}/K_\mathfrak{p}$ is defined for $L_\mathfrak{P} \cap K_\mathfrak{p}^{ab}$, i.e on $D(\mathfrak{P}/\mathfrak{p})/[D(\mathfrak{P}/\mathfrak{p}), D(\mathfrak{P}/\mathfrak{p})]$. Moreover, by the restriction property of higher ramification groups, we have
$$D^{\nu_\mathfrak{p}}(L_\mathfrak{P} \cap K_\mathfrak{p}^{ab}/K_\mathfrak{p}) = \frac{D^{\nu_\mathfrak{p}}(\mathfrak{P}/\mathfrak{p})}{D^{\nu_\mathfrak{p}}(\mathfrak{P}/\mathfrak{p}) \cap [D(\mathfrak{P}/\mathfrak{p}), D(\mathfrak{P}/\mathfrak{p})]}.$$

Also, $\rho_\mathfrak{p}(U_\mathfrak{p}^{(\nu_\mathfrak{p})}) = D^{\nu_\mathfrak{p}}(L_\mathfrak{P} \cap K_\mathfrak{p}^{ab}/K_\mathfrak{p})$; this gives the correspondence between the natural filtration of units and the higher ramification groups.

Since $D(\mathfrak{P}/\mathfrak{p})/D^0(\mathfrak{P}/\mathfrak{p})$ is cyclic, we have
$$[D(\mathfrak{P}/\mathfrak{p}), D(\mathfrak{P}/\mathfrak{p})] = [D(\mathfrak{P}/\mathfrak{p}), D^0(\mathfrak{P}/\mathfrak{p})].$$

This seemingly minor simplification will turn out to be quite important, as we will see in a moment. Recall that $K_S$ is the maximal $p$-extension of $K$ unramified outside $S$, with Galois group $G_S = \text{Gal}(K_S/K)$; it will be convenient here to think of $K_S$ as $K_{S,\nu}$ with $\nu_\mathfrak{p} = \infty$ for $\mathfrak{p} \in S$. We extend $\nu$ to all places of $K$ by $\nu_\mathfrak{p} = 0$ for $\mathfrak{p} \notin S$.

Let $\mathcal{G}$ be the Galois group over $K$ of its maximal $p$-extension $\hat{K}$. We let $D_\mathfrak{p}$ be the Galois group of the maximal $p$-extension $\hat{K}_\mathfrak{p}$ of $K_\mathfrak{p}$. Let $N_{S,\nu}$ be the closed normal subgroup of $\mathcal{G}$ generated by all higher ramification groups $D_\mathfrak{p}^{\nu_\mathfrak{p}}$, as $\mathfrak{p}$ runs over all places of $K$.

**Definition 5.3.** Let $B_{S,\nu} = \ker\left(H^2(G_{S,\nu}, \mathbb{F}_p) \to H^2(\mathcal{G}, \mathbb{F}_p)\right).$

The exact sequence
$$1 \longrightarrow N_{S,\nu} \longrightarrow \mathcal{G} \longrightarrow G_{S,\nu} \longrightarrow 1$$

gives, via Hochschild-Serre, the exact sequence

$$1 \longrightarrow B_{S,\nu}^* \longrightarrow \frac{N_{S,\nu}}{N_{S,\nu}^p[\mathcal{G}, N_{S,\nu}]} \xrightarrow{\alpha} \frac{\mathcal{G}}{\mathcal{G}^p[\mathcal{G}, \mathcal{G}]},$$

allowing us to identify $\ker(\alpha)$ with $B_{S,\nu}^*$, the dual of $B_{S,\nu}$.

**Proposition 5.4.** *1) There is a natural injection*

$$H^2(\mathcal{G}, \mathbb{F}_p) \hookrightarrow \oplus_{\mathfrak{p}} H^2(D_{\mathfrak{p}}, \mathbb{F}_p),$$

*where the sum extends over all primes $\mathfrak{p}$ of $K$. Moreover if $K$ contains $\zeta_p$, we can omit an arbitrary prime $\mathfrak{p}$ in the sum.*

*2) The group $D_{\mathfrak{p}}$ is free (one-relator) if $K_{\mathfrak{p}}$ does not (does) contain $\zeta_p$. In other words, $p\text{-rk} H^2(D_{\mathfrak{p}}, \mathbb{F}_p) = \delta_{\mathfrak{p}}$.*

*3) In the composite map*

$$H^2(G_{S,\nu}, \mathbb{F}_p) \longrightarrow H^2(\mathcal{G}, \mathbb{F}_p) \longrightarrow \oplus_{\mathfrak{p}} H^2(D_{\mathfrak{p}}, \mathbb{F}_p)$$

*the image of $H^2(G_{S,\nu}, \mathbb{F}_p)$ in $\oplus_{\mathfrak{p} \notin S} H^2(D_{\mathfrak{p}}, \mathbb{F}_p)$ is trivial.*

*Proof.* This is proved in [K2]. For part 3), Koch treats only $G_S$ but the same argument applies to $G_{S,\nu}$. See also [Se2], [Sh], [K1]. □

We have therefore proved

**Theorem 5.5.** *For a finite set $S$ of places of $K$ and an arbitrary indexing function $\nu$ of $S$, we have*

$$r_{S,\nu} = p\text{-rk} H^2(G_{S,\nu}, \mathbb{F}_p) \leq p\text{-rk} B_{S,\nu} + \theta_S - \delta_K + \sum_{\mathfrak{p} \in S} \delta_{\mathfrak{p}}.$$

Thus, in order to complete the proof of theorem 5.2, (or to answer Question 5.1), it suffices to bound $p\text{-rk} B_{S,\nu}$ from above. Let us attempt to follow the proof given by Shafarevich and Koch for bounding $B_S$ by relating it to the kernel of the natural map

$$\eta : \mathcal{U}_{S,\nu}/\mathcal{U}_{S,\nu}^p \longrightarrow J_K/K^\times J_K^p.$$

Let us observe, first of all, that $\ker(\eta)$ is under control for arbitrary $\nu$.

**Lemma 5.6.** *We have*

$$p\text{-rk}\ker(\eta) \leq p\text{-rk} E_K + d_{S,\nu} + \sum_{\substack{\mathfrak{p} \in S \cap S_p \\ \nu_{\mathfrak{p}} \neq \infty}} \{[K_{\mathfrak{p}} : \mathbb{Q}_p] + \delta(\mathfrak{p}, \nu_{\mathfrak{p}})\} - \sum_{\mathfrak{p} \in S - S_p} \delta_{\mathfrak{p}} - \sum_{\mathfrak{p} \in S \cap S_p} p\text{-rk} \frac{U_{\mathfrak{p}}^{(1)}}{U_{\mathfrak{p}}^{(\nu_{\mathfrak{p}})}},$$

*where $\delta(\mathfrak{p}, \nu_{\mathfrak{p}})$ is 1 if $\zeta_p \in U_{\mathfrak{p}}^{(\nu_{\mathfrak{p}})}$ and is 0 otherwise.*

*Proof.* There is an exact sequence

$$1 \longrightarrow \frac{\mathcal{U}_{S,\nu} \cap \mathcal{U}_{\emptyset}^p}{\mathcal{U}_{S,\nu}^p} \longrightarrow \ker(\eta) \xrightarrow{\phi} \Delta_{S,\nu} \longrightarrow 1,$$

where $\phi$ is defined as follows. We have

$$\ker(\eta) = \frac{\mathcal{U}_{S,\nu} \cap K^\times J_K^p}{\mathcal{U}_{S,\nu}^p}.$$

Let $u \in \mathcal{U}_{S,\nu} \cap K^\times J_K^p$; write $u = x.j^p$ with $x \in K^\times$, $j \in J_K$. We let $\phi(u)$ be the image of $x$ in $\Delta_{S,\nu}$. The kernel of $\phi$ is easily seen to be $(\mathcal{U}_{S,\nu} \cap \mathcal{U}_\emptyset^p)/\mathcal{U}_{S,\nu}^p$, which is trivial when $\nu(S) \subseteq \{0, \infty\}$, and whose $p$-rank in general is bounded by

$$\begin{aligned} p\text{-rk}\ker(\phi) &= \sum_{\mathfrak{p} \in S \cap S_p} p\text{-rk}(U_\mathfrak{p}^{(\nu_\mathfrak{p})} \cap U_\mathfrak{p}^p) \\ &\leq \sum_{\substack{\mathfrak{p} \in S \cap S_p \\ \nu_\mathfrak{p} \neq \infty}} p\text{-rk}U_\mathfrak{p}^{(\nu_\mathfrak{p})} = \sum_{\substack{\mathfrak{p} \in S \cap S_p \\ \nu_\mathfrak{p} \neq \infty}} \{[K_\mathfrak{p} : \mathbb{Q}_p] + \delta(\mathfrak{p}, \nu_\mathfrak{p})\}. \end{aligned}$$

We now only need apply the formula of Theorem 3.7. □

There are natural surjections

$$\mathcal{U}_{S,\nu} \xrightarrow{\prod_\mathfrak{p} \rho_\mathfrak{p}} \prod_\mathfrak{p} \frac{D_\mathfrak{p}^{\nu_\mathfrak{p}}}{[D_\mathfrak{p}, D_\mathfrak{p}] \cap D_\mathfrak{p}^{\nu_\mathfrak{p}}} \xleftarrow{\beta} \prod_\mathfrak{p} \frac{D_\mathfrak{p}^{\nu_\mathfrak{p}}}{[D_\mathfrak{p}, D_\mathfrak{p}^{\nu_\mathfrak{p}}]} \xrightarrow{\lambda} \frac{N_{S,\nu}}{[\mathcal{G}, N_{S,\nu}]}. \tag{1}$$

Note that $\beta$ is "pointing the wrong way." In the classical case (where $\nu_\mathfrak{p}$ is either 0 or $\infty$), $\beta$ is an isomorphism ($[D_\mathfrak{p}, D_\mathfrak{p}] = [D_\mathfrak{p}, D_\mathfrak{p}^0]$), so we can define a surjective global map

$$\mathcal{U}_S \xrightarrow{\rho} N_S/N_S^p[\mathcal{G}, N_S].$$

We then have a commutative diagram:

$$\begin{array}{ccccc} B_S^* & \hookrightarrow & \dfrac{N_S}{N_S^p[\mathcal{G}, N_S]} & \xrightarrow{\alpha} & \dfrac{\mathcal{G}}{\mathcal{G}^p[\mathcal{G}, \mathcal{G}]} \\ {\scriptstyle \exists \rho''} \Big\uparrow & & {\scriptstyle \rho}\Big\uparrow & & {\scriptstyle \rho'}\Big\uparrow \\ \ker(\eta) & \hookrightarrow & \dfrac{\mathcal{U}_S}{\mathcal{U}_S^p} & \xrightarrow{\eta} & \dfrac{J_K}{K^\times J_K^p} \end{array}$$

where $\rho$ is as defined above. As $\rho'$ is injective, there is a natural surjection $\rho'' : \ker(\eta) \longrightarrow B_S^*$. Putting this together with Theorem 5.5 and Lemma 5.6, we have proved Theorem 5.2.

To summarize, in the general case, what is lacking is a bound for the dimension of $B_{S,\nu}$, for which it would suffice to have a natural map $\rho^? : \mathcal{U}_{S,\nu}/\mathcal{U}_{S,\nu}^p \longrightarrow N_{S,\nu}/N_{S,\nu}^p[\mathcal{G}, N_{S,\nu}]$ with finite cokernel. The classical approach sketched above does not immediately generalize because the map $\beta$ (an isomorphism in the classical case) has infinite kernel in general, by the following remarkable theorem of Gordeev [G].

**Theorem 5.7** (Gordeev). *Let $K_\mathfrak{p}$ be a finite extension of $\mathbb{Q}_p$. Let $\hat{K}_\mathfrak{p}$ be the maximal $p$-extension of $K_\mathfrak{p}$ and put $D_\mathfrak{p} = \mathrm{Gal}(\hat{K}_\mathfrak{p}/K_\mathfrak{p})$. For $\nu_\mathfrak{p} > 1$,*

$$r(D_\mathfrak{p}/D_\mathfrak{p}^{\nu_\mathfrak{p}}) = \dim_{\mathbb{F}_p} H^2(D_\mathfrak{p}/D_\mathfrak{p}^{\nu_\mathfrak{p}}, \mathbb{F}_p) = \infty.$$

**Corollary 5.8.** *For each $K$-prime $\mathfrak{p}$, put*

$$X_\mathfrak{p} = \frac{D_\mathfrak{p}^{\nu_\mathfrak{p}}}{[D_\mathfrak{p}, D_\mathfrak{p}^{\nu_\mathfrak{p}}]}, \qquad Y_\mathfrak{p} = \frac{D_\mathfrak{p}^{\nu_\mathfrak{p}}}{[D_\mathfrak{p}, D_\mathfrak{p}] \cap D_\mathfrak{p}^{\nu_\mathfrak{p}}}.$$

*Let $\beta : \prod_\mathfrak{p} X_\mathfrak{p} \longrightarrow \prod_\mathfrak{p} Y_\mathfrak{p}$ be the natural map appearing in (1). Suppose $\nu_\mathfrak{p} > 1$ for some $\mathfrak{p} \in S \cap S_p$. Then $\ker(\beta)$ is infinite. Indeed, with hypotheses as in the theorem,*

$$p\text{-rk}X_\mathfrak{p} = \infty, \qquad p\text{-rk}Y_\mathfrak{p} \leq [K_\mathfrak{p} : \mathbb{Q}_p] + 1.$$

*Proof.* Let us take the exact sequence
$$1 \longrightarrow D_{\mathfrak{p}}^{\nu_{\mathfrak{p}}} \longrightarrow D_{\mathfrak{p}} \longrightarrow D_{\mathfrak{p}}/D_{\mathfrak{p}}^{\nu_{\mathfrak{p}}} \longrightarrow 1.$$

Applying the Hochschild-Serre sequence, we get an injection
$$\left(\ker\left(H^2(D_{\mathfrak{p}}/D_{\mathfrak{p}}^{\nu_{\mathfrak{p}}}, \mathbb{F}_p) \to H^2(D_{\mathfrak{p}}, \mathbb{F}_p))\right)\right)^* \hookrightarrow \frac{D_{\mathfrak{p}}^{\nu_{\mathfrak{p}}}}{[D_{\mathfrak{p}}, D_{\mathfrak{p}}^{\nu_{\mathfrak{p}}}](D_{\mathfrak{p}}^{\nu_{\mathfrak{p}}})^p}.$$

But $D_{\mathfrak{p}}$ is a $p$-group with $p$-rk$H^2(D_{\mathfrak{p}}, \mathbb{F}_p) = 0$ or $1$, whereas by Theorem 5.7, $D_{\mathfrak{p}}/D_{\mathfrak{p}}^{\nu_{\mathfrak{p}}}$ has infinite relation-rank. Thus, the $p$-rank of $\frac{D_{\mathfrak{p}}^{\nu_{\mathfrak{p}}}}{[D_{\mathfrak{p}}, D_{\mathfrak{p}}^{\nu_{\mathfrak{p}}}](D_{\mathfrak{p}}^{\nu_{\mathfrak{p}}})^p}$ is infinite. On the other hand, the quotient $\frac{D_{\mathfrak{p}}^{\nu_{\mathfrak{p}}}}{[D_{\mathfrak{p}}, D_{\mathfrak{p}}] \cap D_{\mathfrak{p}}^{\nu_{\mathfrak{p}}}(D_{\mathfrak{p}}^{\nu_{\mathfrak{p}}})^p}$ is the restriction of $D_{\mathfrak{p}}^{\nu_{\mathfrak{p}}}$ to the maximal elementary abelian $p$-extension of $K_{\mathfrak{p}}$, and so has $p$-rank at most $[K_{\mathfrak{p}} : \mathbb{Q}_p] + 1$.  □

*Remark* 5.9. Note that this approach can be salvaged and a relation-rank bound can be obtained if $\ker(\beta)$ and $\ker(\lambda)$ are commensurable. For instance, when do we have $\ker(\beta) \subseteq \ker(\lambda)$?

### 5.2. Not all $G_{S,\nu}$ are finitely presentable.

The groups $G_{S,\nu}$ are finitely generated, so they are finitely presentable in the category of pro-$p$ groups if and only if they have finite Euler characteristic $\chi_2(G_{S,\nu}) = r_{S,\nu} - d_{S,\nu}$. Note that, for *finitely indexed* $(S, \nu)$, $\chi_2(G_{S,\nu}) \geq 0$ since $G_{S,\nu}$ has finite abelianization in that case.

The main result of this part gives infinitely many examples (with $p = 2$, $\nu_{\mathfrak{p}}$ finite for some $\mathfrak{p}$ but infinite for others) for which $r_{S,\nu} = \chi_2(G_{S,\nu}) = \infty$.

**Theorem 5.10.** *Suppose $p = 2$, and $\ell$ is a prime with $\ell \equiv 7 \pmod{16}$. Put $K = \mathbb{Q}(\sqrt{-\ell})$, in which $2$ splits: $2\mathcal{O}_K = \mathfrak{p}\mathfrak{q}$. Choose $i > 1$, and define $\nu : S = \{\mathfrak{p}, \mathfrak{q}\} \to [0, \infty]$ by $\nu_{\mathfrak{p}} = i, \nu_{\mathfrak{q}} = \infty$. Then $r(G_{S,\nu}) = \infty$.*

*Remark* 5.11. 1) If we take $\ell = 7$ and $\nu_{\mathfrak{p}} = \nu_{\mathfrak{q}} = 2$, then $G_{S,\nu}^{ab}$ has order $2$ hence so does $G_{S,\nu}$; that $G_{S,\nu}$ is finite in this case can be deduced from Theorem 4.2 and discriminant bounds as well [O]. Assuming the Generalized Riemann Hypothesis, the latter method also gives the finiteness of $G_{S,\nu}$ for $\nu_{\mathfrak{p}} = \nu_{\mathfrak{q}} = 3$ (here $G_{S,\nu}^{ab}$ is $(2, 2, 2)$). For $\nu_{\mathfrak{p}} = \nu_{\mathfrak{q}} = 4$, $G_{S,\nu}^{ab}$ is $(4, 4, 2)$, and we do not know whether $G_{S,\nu}$ is infinite or not.

2) Note that $G_{S,\nu}$ (with $(S, \nu)$ as in the theorem) has $G_{\{\mathfrak{q}\}}$ has as a quotient, which in turn has $\mathbb{Z}_2$ as a quotient.

To prove Theorem 5.10 we have to recall a result of Wingberg [W].

**Definition 5.12.** Recall that $D_{\mathfrak{p}} = \mathrm{Gal}(\hat{K}_{\mathfrak{p}}/K_{\mathfrak{p}})$ where $\hat{K}_{\mathfrak{p}}$ is the maximal $p$-extension of $K_{\mathfrak{p}}$. For a number field $K$, a $p$-extension $L/K$ is called *local* (*maximal local*) at $\mathfrak{p}$ if the natural composite map
$$D_{\mathfrak{p}} \xrightarrow{\mathrm{res}} \mathrm{Gal}(L_{\mathfrak{P}}/K_{\mathfrak{p}}) \hookrightarrow \mathrm{Gal}(L/K)$$
is surjective (an isomorphism).

Note that if $L/K$ is of local, then for every tower $L/F/K$, where $F/K$ is Galois, $F/K$ is local. We introduce a little more notation. For a finite set $S$ of primes of $K$, let $\overline{S}$ be $S$ if $p > 2$ and the union of $S$ with the real infinite places of $K$ if $p = 2$. Then $K_{\overline{S}}/K$ is the maximal $p$-extension of $K$ unramified outside $\overline{S}$ (the real infinite places are allowed to ramify in $K_{\overline{S}}$ when $p = 2$). Under the condition $S_p \subseteq S$, Wingberg characterizes those $K_{\overline{S}}/K$ which are of local or maximal local type ([W], Corollary 1.5 and Theorem 1.6).

**Theorem 5.13** (Wingberg). *Let $S, \overline{S}$ be as above. Let $r_2$ denote the number of imaginary places of $K$. Then, $K_{\overline{S}}/K$ is local at $\mathfrak{p} \in S_p$ if and only if*

$$\sum_{\mathfrak{q} \in S - \{\mathfrak{p}\}} \delta_\mathfrak{q} - \delta_K + p\text{-rk}\Delta^{\{\mathfrak{p}\}}_{S-\{\mathfrak{p}\}} + r_2 = [K_\mathfrak{p} : \mathbb{Q}_p].$$

*Moreover, $K_{\overline{S}}/K$ is maximal local at $\mathfrak{p} \in S_p$ if and only if the following conditions are satisfied:*

(1) *If $p = 2$, $K$ is totally imaginary,*
(2) $\sum_{\mathfrak{q} \in S \setminus \{\mathfrak{p}\}} \delta_\mathfrak{q} = \delta_K$,
(3) $\Delta^{\{\mathfrak{p}\}}_{S-\{\mathfrak{p}\}} = \{1\}$,
(4) $r_2 = [K_\mathfrak{p} : \mathbb{Q}_p]$.

**Corollary 5.14.** *With notation as above, and $p = 2$, $K_{\overline{S}}/K$ is maximal local at $\mathfrak{p} \in S_p$ if and only if*

(1) *$K$ is totally imaginary,*
(2) *$S = \{\mathfrak{p}, \mathfrak{q}\}$ where $2\mathcal{O}_K = \mathfrak{p}\mathfrak{q}$, and $\mathfrak{p} \neq \mathfrak{q}$,*
(3) *There does not exist a quadratic extension of $K$ unramified outside $\mathfrak{p}$ in which $\mathfrak{q}$ splits,*
(4) $r_2 = [K_\mathfrak{p} : \mathbb{Q}_p]$.

*Proof.* Since $-1 = \zeta_2 \in K$, $S - \{\mathfrak{p}\}$ must be a singleton. Moreover, by Proposition 1.1, a quadratic extension of $K$, $\mathfrak{q}$-decomposed and unramified outside $\mathfrak{p}$, exists if and only if $\Delta^{\{\mathfrak{p}\}}_{S-\{\mathfrak{p}\}} \neq \{1\}$. □

We need to verify condition 3 of the above corollary for the fields appearing in Theorem 5.10. We do so in the following Proposition, which can be proved directly via a long calculation; we will instead derive it more easily by making use of a result in a forthcoming book of Gras [Gr].

**Proposition 5.15.** *Let $K = \mathbb{Q}(\sqrt{-\ell})$ where $\ell \equiv 7 \mod 8$ is a prime, $2\mathcal{O}_K = \mathfrak{p}\mathfrak{q}$. There exists a quadratic extension of $K$, unramified outside $\mathfrak{p}$, in which $\mathfrak{q}$ splits completely if and only if $\ell \equiv 15 \mod 16$.*

*Proof.* Let $n$ be the order of $\mathfrak{p}$ in $\text{Cl}_K$: by genus theory, $n$ is odd. A generator $\alpha \in \mathcal{O}_K$ of $\mathfrak{p}^n$ is of the form $\alpha = (a + b\sqrt{-\ell})/2$, where $a$ and $b$ are odd integers. Call the non-trivial element of $\text{Gal}(K/\mathbb{Q})$ $\tau$, so that $\alpha^\tau = (a - b\sqrt{-\ell})/2$, and $\mathfrak{q}^n = (\alpha^\tau)$. The case $\ell = 7$ is easily checked by hand. Now suppose $\ell > 7$; then 2 is not a norm in $K/\mathbb{Q}$, so $n \geq 3$. We find a sequence of equivalent conditions as follows.

**Claim.** *For $\ell > 7$, the following are equivalent:*

(1) *There exists a quadratic extension of $K$, unramified outside $\mathfrak{p}$, in which $\mathfrak{q}$ splits completely,*
(2) $[K_\mathfrak{p}(\sqrt{-1}, \sqrt{\alpha^\tau}) : K_\mathfrak{p}] = 2$,
(3) $(2, \alpha^\tau)_\mathfrak{p} = 1$,
(4) $\alpha^\tau \equiv \pm 1 \mod \mathfrak{p}^3$,
(5) $a^2 \equiv 1 \mod 16$,
(6) $\ell \equiv 15 \mod 16$.

*Proof of Claim.* The equivalence of 1 and 2 is from Gras' book (Chapter V). We then have $\alpha^\tau \in K_\mathfrak{p}(\sqrt{-1}) \iff \alpha^\tau \equiv \pm 1 \pmod{\mathfrak{p}^3} \iff (2, \alpha^\tau)_\mathfrak{p} = 1$, because $\alpha^\tau$ is a local unit at $\mathfrak{p}$ (see [Se1] for example). Since $n \geq 3$, condition 4 is equivalent to $\alpha \equiv 0 \pmod{\mathfrak{p}^3}$ and $\alpha + \alpha^\tau \equiv \pm 1 \pmod{\mathfrak{p}^3}$, i.e. to $a \equiv \pm 1 \pmod 8$ or $a^2 \equiv 1 \bmod 16$. For the final step, we want to show that this is equivalent to $\ell \equiv 15 \bmod 16$. We take the norm of $2\alpha$: $4\mathbb{N}\alpha = a^2 + b^2\ell = 2^{n+2}$, to see that $2^{n+2}$ is a square mod $b$, and as $n$ is odd, 2 is a square mod $b$. Now by the quadratic reciprocity law, 2 is a square mod $b$ implies $b \equiv \pm 1 \bmod 8$. We again use the fact that $n \geq 3$: $a^2 + b^2\ell = 2^{n+2}$ gives $a^2 + b^2\ell \equiv 0 \pmod{16}$ or $a^2 + \ell \equiv 0 \pmod{16}$. This concludes the proof of the claim, and of Proposition 5.15. □

*Proof.* [Theorem 5.10] With $S = \{\mathfrak{p}, \mathfrak{q}\}$, $K_S/K$ is maximal local at $\mathfrak{q}$ by Proposition 5.15 and Corollary 5.14. Thus, the global Galois group $G_{S,\nu}$ is isomorphic to the local group appearing in Gordeev's Theorem 5.7, which has infinite relation rank. □

5.3. **Variation of $r_{S,\nu}$ for fixed $S$.** Let us consider what happens as we vary $\nu$ for $K = \mathbb{Q}(\sqrt{-\ell})$, $\ell \equiv 7 \bmod 16$, $p = 2$, $S = S_p = \{\mathfrak{p}, \mathfrak{q}\}$ as above. We introduce some notation. For integers $i, j > 1$, let $\nu_{i,j} : S = \{\mathfrak{p}, \mathfrak{q}\} \to [0, \infty)$ be defined by $\nu_{i,j}(\mathfrak{p}) = i$, $\nu_{i,j}(\mathfrak{q}) = j$. To simplify, we write

$$G_{i,j} = G_{S,\nu_{i,j}}.$$

**Corollary 5.16.** *For fixed $i$, we have $\sup_j r(G_{i,j}) = \sup_j \chi_2(G_{i,j}) = \infty$.*

*Proof.* Fix $i$. We have $G_{i,\infty} = \varprojlim_j G_{i,j}$. Consequently,

$$H^2(G_{i,\infty}, \mathbb{F}_p) = \varinjlim_j H^2(G_{i,j}, \mathbb{F}_p).$$

If for some $j$, $r(G_{i,j}) = \infty$, then there is nothing to prove. Now suppose $r(G_{i,j}) \leq M < \infty$ for all $j$. Then $H^2(G_{i,\infty}, \mathbb{F}_p)$ is a direct limit of abelian groups of cardinality at most $M$, hence is finite. This contradicts Theorem 5.10 ($r(G_{i,\infty}) = \infty$). The claim for $\sup_j \chi_2(G_{i,j})$ follows since $p\text{-rk} G_{i,j} \leq p\text{-rk} G_{i,\infty} \leq p\text{-rk} G_S < \infty$. □

*Remark* 5.17. Note that $\sup_j(G_{i,j}) = \infty$ implies that $G_{i,\infty}$ is not $p$-adic analytic. Consider $G_S$, the Galois group over $K = \mathbb{Q}(\sqrt{-\ell})$ (with prime $\ell \equiv 7 \bmod 16$)) of the maximal 2-extension unramified outside 2. Then $G_S$ is not analytic, since it has non-analytic quotients $G_{i,\infty}$. Note that we have $G_S = \varprojlim_i G_{i,\infty}$, with $\chi_2(G_{i,\infty}) = \infty$, but $\chi_2(G_S) = -(r_2 + 1) = -2$ is finite.

Part 3. **The Fontaine-Mazur conjecture and $p$-adic representations of $G_{S,\nu}$**

6. THEOREMS OF SEN AND COATES-GREENBERG

Suppose $F$ is a complete local field of characteristic zero with valuation $v_F$ and residue field of characteristic $p$. Let $E/F$ be a totally ramified Galois extension such that $G = \text{Gal}(E/F)$ is an infinite $p$-adic analytic group. The group $G$ has two natural filtrations given by its $p$-central series ($G_1 = G$, $G_{j+1} = G_i^p[G, G_j]$), and the higher ramification groups $D^i(E/F)$. It was conjectured by Serre and proven by Sen [S1] that these filtrations are closely related as follows.

**Theorem 6.1** (Sen)**.** *Let $F$ be a complete local field of characteristic zero with residue field of characteristic $p$ and ramification index $e = v_F(p)$. Suppose $E/F$ is a totally ramified $p$-extension with $p$-adic analytic Galois group $G$. Then there exists a constant $c$, such that for all $n \geq 0$,*

$$D^{ne+c}(E/F) \subseteq G_n \subseteq D^{ne-c}(E/F).$$

The following corollary appears in the paper [CG] of Coates and Greenberg (Theorem 2.13) and was pointed out to us by Schmidt [Sch2]. In an earlier version of this paper, we used a slightly weaker statement along these lines to prove a weaker version of Theorem 7.3. We reproduce the proof, since it is not long.

**Corollary 6.2** (Coates-Greenberg)**.** *Suppose $F_0$ is a finite extension of $\mathbb{Q}_p$ and $E/F_0$ is a $p$-extension with $p$-adic analytic Galois group. Let $F$ be the maximal unramified $p$-extension of $F_0$ contained in $E/F_0$. Assume that the ramification in $E/F_0$ is of bounded depth, i.e. there exists $i > 0$ such that $D^i(E/F_0) = \{1\}$. Then $E/F$ is finite. In other words, the inertia group of $E/F_0$ is finite.*

*Proof.* The group $G = \mathrm{Gal}(E/F)$ is a closed subgroup of a $p$-adic analytic group, hence $p$-adic analytic. First note that $F$ is a characteristic zero local field with a discrete valuation $v$ and a (possibly infinite) perfect residue field of characteristic $p$. Let $\hat{F}$ be the completion of $F$ at $v$. By restriction, $\mathrm{Gal}(E\hat{F}/\hat{F})$ (resp. $D^i(E\hat{F}/\hat{F})$) is isomorphic to $\mathrm{Gal}(E/F)$ (resp. $D^i(E/F)$). As $F/F_0$ is unramified, the groups $D^i(E/F_0)$ and $D^i(E/F)$ are the same, and by assumption, they vanish. By Theorem 6.1 applied to the totally ramified extension $E\hat{F}/\hat{F}$, there exists a constant $c$ such that

$$D^{ne+c}(E\hat{F}/\hat{F}) \subseteq \mathrm{Gal}(E\hat{F}/\hat{F})_n \subseteq D^{ne-c}(E\hat{F}/\hat{F})).$$

Put $m = \lceil (i+c)/e \rceil$ and $j = me - c$. Then $j \geq i$, so $D^j(E\hat{F}/\hat{F})$ is trivial, hence $G_m$ (the $m$th term of the central $p$-series of $G \simeq \mathrm{Gal}(E\hat{F}/\hat{F})$) is trivial as well. Thus, $\mathrm{Gal}(E/F)$ is finite. $\square$

## 7. The Fontaine-Mazur conjecture

Let us call a finitely generated pro-$p$ group $G$ *$p$-adically finite* if it has no infinite $p$-adic analytic quotients. In other words, all $p$-adic representations of $G$ into $\mathrm{GL}_n(\mathbb{Q}_p)$ have finite image. Recall the Fontaine-Mazur conjecture.

**Conjecture 7.1** ([FM], Conjecture 5a)**.** *Let $K$ be a number field, $S$ a finite set of places of $K$ disjoint from $S_p$. Then $G_S$ is $p$-adically finite.*

For partial corroboration of this conjecture, see Boston [B1], [B2], [B3] and Hajir [Ha]. Conjecture 7.1 comes about as a consequence of Fontaine and Mazur's vast program for characterizing the $p$-adic Galois representations which "come from algebraic geometry," meaning those isomorphic to a subquotient of the action of $\mathrm{Gal}(\overline{K}/K)$ on an étale cohomology group $H^q_{\mathrm{ét}}(X_{\overline{K}}, \mathbb{Q}_p(r))$, where $X$ is a smooth projective variety over $K$. They call an irreducible $p$-adic representation of $\mathrm{Gal}(\overline{K}/K)$ "geometric" if is satisfies two conditions: 1) it is unramified outside a finite set of places of $K$, and 2) its restriction to every decomposition group is potentially semistable. Their main conjecture, then, is that an irreducible $p$-adic representation is geometric if and only if it comes from algebraic geometry. One direction of this conjecture ("algebro-geometric representations are potentially semistable") has a much longer history, and is now established (see, for example, Tsuji [T]).

By a theorem of Grothendieck (cf. appendix of [ST]), a tamely ramified $p$-adic representation is always potentially semistable. Thus, when $S \cap S_p = \emptyset$, every $p$-adic representation $\rho$ of $G_S$ should come from algebraic geometry; algebro-geometric considerations (see [FM] and [S2]) then imply that $\rho$ has finite image. This is equivalent to the tamely ramified Fontaine-Mazur Conjecture, i.e. Conjecture 7.1 above.

Following the philosophy outlined in the introduction, we can formulate the following extension of that conjecture to the case of wild ramification of finite depth.

**Conjecture 7.2.** If $(S, \nu)$ is a finitely indexed set, then $G_{S,\nu}$ is $p$-adically finite.

As was pointed out to us by Schmidt, the Coates-Greenberg corollary to Sen's theorem has the following consequence.

**Theorem 7.3.** *Suppose $(S, \nu)$ is a finitely indexed set for $K$. Let $L/K$ be a Galois subextension of $K_{S,\nu}/K$ such that the Galois group $\mathrm{Gal}(L/K)$ is $p$-adic analytic. Then $L/K$ is potentially tamely ramified.*

*Proof.* By Corollary 6.2, for all places $\mathfrak{P}$ of $L$ dividing $\mathfrak{p} \in S \cap S_p$, the inertia group $D^0(L/K, \mathfrak{P})$ is finite. Hence $L/K$ is potentially tamely ramified: there exists a number field $K'$ in $L/K$ such that $L/K'$ is unramified at all places above $p$. $\square$

We have three immediate corollaries.

**Corollary 7.4.** *If $\nu$ is finite, every $p$-adic representation of $\mathrm{Gal}(\overline{K}/K)$ factoring through $G_{S,\nu}$ is potentially semistable.*

*Proof.* Via a finite base change, we pass to a tamely ramified extension and apply Grothendieck's theorem. $\square$

As a consequence, $p$-adic representations of $\mathrm{Gal}(\overline{K}/K)$ factoring through $G_{S,\nu}$ (when $\nu$ is finite) should come from algebraic geometry.

**Corollary 7.5.** *Conjecture 7.1 implies Conjecture 7.2.*

*Proof.* Let $L/K$ be a Galois subextension of $K_{S,\nu}/K$ such that the Galois group $\mathrm{Gal}(L/K)$ is $p$-adic analytic. There exists a number field $K'$ in $L/K$ such that $L/K'$ is unramified at all places above $p$. Moreover $\mathrm{Gal}(L/K')$ is $p$-adic analytic and so Conjecture 7.1 implies the finiteness of $L/K'$, and thus of $L/K$. $\square$

**Corollary 7.6.** *Assume Conjecture 7.1. If $L/K$ is an infinite Galois extension with $p$-adic analytic Galois group $\mathrm{Gal}(L/K)$, then either infinitely many primes of $K$ ramify in $L$, or $L/K$ is deeply ramified at some prime $\mathfrak{p} \in S_\mathfrak{p}$.*

*Proof.* If $L/K$ were ramified at only a finite set of places $S$ of $K$, and the ramification were of bounded depth, then $L/K$ is finite by the previous Corollary. $\square$

Thus, the Fontaine-Mazur conjecture implies that infinite $p$-adic analytic extensions of number fields are infinitely ramified, either horizontally or vertically, so to speak. Only recently have $p$-adic Lie extensions with infinitely many ramified primes been constructed: see the work of Ramakrishna, Khare, and Rajan [R], [KRa], [KR]. Recalling Theorem 4.4, we can reformulate The Fontaine-Mazur Conjecture as follows.

**Corollary 7.7.** *Conjecture 7.1 is equivalent to the following statement: An infinite $p$-adic analytic extension of a number field is asymptotically bad.*

EXTENSIONS OF NUMBER FIELDS WITH WILD RAMIFICATION OF BOUNDED DEPTH   21*Proof.* Suppose $K$ is a number field and $L/K$ is an infinite $p$-adic analytic extension; if we admit Fontaine-Mazur, then, by Corollary 7.6, $L/K$ is either ramified at infinitely many primes (in which case it is easy to see that it is asymptotically bad) or it is deeply ramified (in which case it is asymptotically bad by Theorem 4.4). On the other hand, supposing that every $p$-adic analytic extension of $K$ is asymptotically bad, and knowing (Theorem 4.2) that an infinite tame extension unramified outside a finite set of primes is asymptotically good, we would conclude that $K$ admits no infinite tame extension unramified outside a finite set of primes. □

In the remainder of this section, we give a criterion for $(S, \nu)$ to satisfy Conjecture 7.2, and provide two kinds of unconditional examples where this criterion, and therefore Conjecture 7.2, hold. The main tool, once again, is Wingberg's study of local and global groups.

**Theorem 7.8.** *If $K_{S,\nu}/K$ is of local type at some $\mathfrak{p} \in S \cap S_p$, and $\nu$ is finite, then $G_{S,\nu}$ is $p$-adically finite.*

*Proof.* Consider a Galois extension $F$ of $K$ contained in $K_{S,\nu}$, such that $\mathrm{Gal}(F/K)$ is $p$-adic analytic. By the remark following Definition 5.12, one has $\mathrm{Gal}(F/K) \simeq \mathrm{Gal}(F_\mathfrak{P}/K_\mathfrak{p})$. By corollary 6.2, the inertia group $D^0(F/K, \mathfrak{P})$ is finite. Letting $M$ be the subfield of $F/K$ fixed by $D^0(F/K, \mathfrak{P})$, the Galois group $\mathrm{Gal}(M/K)$ is isomorphic to $\mathrm{Gal}(F_\mathfrak{P}/K_\mathfrak{p})/D^0(F/K, \mathfrak{P})$, which is abelian. In particular, $M$ is a subfield of $K_{S,\nu}^{ab}$, which is a ray class field of finite conductor, hence of finite degree over $K$. Therefore, $M/K$ is finite, and so is $F/K$. □

Now we will give two immediate applications of Theorem 7.8 and Wingberg's Theorem 5.13, which provide examples where Conjecture 7.2 holds. Note that we do not use Conjecture 7.1 in these examples, but by the same token, we do not know whether the groups $G_{S,\nu}$ in question are infinite; thus, it is not yet clear whether these examples represent non-trivial evidence for Conjecture 7.2.

**Corollary 7.9.** *We consider the situation of Theorem 5.10: $p = 2$, $K = \mathbb{Q}(\sqrt{-\ell})$, $\ell \equiv 7 \bmod 16$ is prime, $2\mathcal{O}_K = \mathfrak{p}\mathfrak{q}$. Take $0 < i, j < \infty$, $S = \{\mathfrak{p}, \mathfrak{q}\}$, and put $\nu_{i,j}(\mathfrak{p}) = i$, $\nu_{i,j}(\mathfrak{q}) = j$. Then $G_{S,\nu_{i,j}}$ is 2-adically finite. (For large $i, j$, $G_{S,\nu_{i,j}}$ is a group with three generators).*

**Corollary 7.10.** *Let $p$ be an odd regular prime, and put $K = \mathbb{Q}(\zeta_p)$. We have the factorization: $p\mathcal{O}_K = \mathfrak{p}^{p-1}$. Let $0 < i < \infty$, $S = \{\mathfrak{p}\}$ and $\nu_i$ defined by $\nu_i(\mathfrak{p}) = i$. Then the group $G_{S,\nu_i}$ is $p$-adically finite. (For large $i$, $G_{S,\nu_i}$ is a group with $(p+1)/2$ generators).*

*Remark* 7.11. The group-theoretical method of Boston [B1] and [B2] generalizes from the tame case to the case of bounded-depth wild ramification to provide further evidence for Conjecture 7.2.

## References

[AM] B. Angles and C. Maire, A note on tamely ramified towers of global function fields, *Finite Fields and Appl.*, to appear.
[B1] N. Boston, Some cases of the Fontaine-Mazur conjecture, J. Number Theory **42** (1992), 285-291.
[B2] N. Boston, Some cases of the Fontaine-Mazur conjecture II, J. Number Theory **75** (1999), 161-169.
[B3] N. Boston, $p$-adic Galois representations and pro-$p$ Galois groups, pp. 329-348 in [SSS].
[CG] J. Coates and R. Greenberg, Kummer theory for abelian varieties over local fields, Invent. Math **124** (1996), 129–174.

Farshid Hajir
Dept. of Mathematics
California State University, San Marcos
San Marcos CA 92096
fhajir@csusm.edu

Christian Maire
Laboratoire A2X
Université Bordeaux I
Cours de la Libération
33405 Talence Cedex France
maire@math.u-bordeaux.fr